\documentclass{article}

\usepackage{amsfonts}
\usepackage{graphicx}
\usepackage[latin1]{inputenc}
\usepackage[english]{babel}
\usepackage{color}
\usepackage{float}
\usepackage{fancyhdr}
\usepackage{fancybox}
\usepackage[english]{minitoc}
\usepackage{amssymb,amsmath}
\usepackage{textcomp}
\usepackage{amsmath}
\usepackage{boxedminipage}
\usepackage{minitoc}
\usepackage{lscape}
\usepackage{amsthm}
\usepackage{parskip}
\usepackage{epstopdf}
\usepackage[labelfont=bf]{caption}
\usepackage{subfigure}
\usepackage[top=90pt, bottom=90pt]{geometry}
\usepackage{adjustbox}
\usepackage{longtable}
\usepackage{booktabs}
\usepackage{colortbl}
\usepackage{multirow}
\usepackage{placeins}
\definecolor{grijs}{RGB}{232,232,232}
\usepackage[table]{xcolor}
\usepackage{tikz}
\usetikzlibrary{matrix,backgrounds}
\usepackage{mathabx} 
\numberwithin{equation}{section}

\newtheorem{theorem}{Theorem}[section]

\theoremstyle{definition}

\newcommand{\half}{\frac{1}{2}}
\newcommand{\dud}[1]{\frac{\partial v}{\partial #1}}
\newcommand{\dudd}[1]{\frac{\partial^2v}{\partial #1^2}}
\newcommand{\duddm}[2]{\frac{\partial^2v}{\partial #1 \partial #2}}
\newcommand{\dd}{\mathrm{d}}
\newcommand{\rhoh}{\widehat{\rho}}

\newcommand{\wV}{\widehat{V}}
\def\-{\raisebox{.75pt}{-}}

\def\R{\mathbb{R}}

\pgfdeclarelayer{myback}
\pgfsetlayers{myback,background,main}

\tikzset{mycolor/.style = {line width=1bp,color=#1}}%
\tikzset{myfillcolor/.style = {draw,fill=#1}}%

\author{Lynn Boen\footnote{Department of Mathematics,
University of Antwerp, Middelheimlaan 1, B-2020 Antwerp, Belgium.
\mbox{Email}: \texttt{\{lynn.boen,karel.inthout\}@uantwerpen.be}.}
~and Karel~J.~in 't Hout\footnotemark[\value{footnote}]
}

\title{Operator splitting schemes for American options\\ 
under the two-asset Merton jump-diffusion model}

\begin{document}
\maketitle

\begin{abstract}
This paper deals with the efficient numerical solution of the two-dimensional partial integro-differential complementarity problem (PIDCP) that holds for the value of American-style options under the two-asset Merton jump-diffusion model.
We consider the adaptation of various operator splitting schemes of both the implicit-explicit (IMEX) and the alternating direction implicit (ADI) kind that have recently been studied for partial integro-differential equations (PIDEs) in \cite{BH19}. 
Each of these schemes conveniently treats the nonlocal integral part in an explicit manner.
Their adaptation to PIDCPs is achieved through a combination with the Ikonen--Toivanen splitting technique \cite{IT04} as well as with the penalty method \cite{ZFV98}.
The convergence behaviour and relative performance of the acquired eight operator splitting methods is investigated in extensive numerical experiments for American put-on-the-min and put-on-the-average options.
\end{abstract}

\section{Introduction}
The flexibility of American-style options, which allows the holder to exercise at any time up to and including maturity, renders them popular financial contracts. 
This flexibility poses a challenge in the numerical valuation of these options via partial differential equations (PDEs), since the early exercise feature leads to a nonlinear free boundary problem.
When the underlying asset price process exhibits jumps, this free boundary problem forms a partial integro-differential complementarity problem (PIDCP) where the integral
part is nonlocal. 
In addition, when there are multiple underlying assets, the PIDCP is multidimensional.
The present paper is concerned with the efficient numerical solution of this advanced type of problems.

Classical methods for solving complementarity problems that arise in the valuation of American-style options are the Brennan--Schwartz algorithm~\cite{BS77} and the projected successive over-relaxation method.
As it turns out, however, these two methods are of limited practical use as the former is only applicable under certain, restrictive conditions and the latter, iterative method converges too slowly.
Over the recent years, a variety of effective numerical methods has been developed in the literature for P(I)DCPs that model American option values under various asset price processes. 
We provide a brief and nonexhaustive overview.

Clarke \& Parrott \cite{CP99} considered the two-dimensional PDCP for American option values under the Heston model and applied a multigrid method.
This approach was next investigated in for example Oosterlee~\cite{O03} and Toivanen \& Oosterlee~\cite{TO12}.

Zvan, Forsyth \& Vetzal~\cite{ZFV98} proposed the penalty method for solving the Heston PDCP. 
The properties of this method were rigorously analyzed in Forsyth \& Vetzal~\cite{FV02} for the one-dimensional Black--Scholes PDCP.
The penalty method was generalized by d'Halluin et al.~\cite{HFL04,HFV05} to one-dimensional PIDCPs for American option values, in particular under the Merton and Kou models, and subsequently applied by Clift \& Forsyth~\cite{CF08} to two-dimensional PIDCPs.
Here a fixed-point iteration is used to efficiently handle the integral part.

Ikonen \& Toivanen~\cite{IT04} introduced an alternative approach for solving the one-dimensional Black--Scholes PDCP.
Here the PDCP is reformulated by means of an auxiliary variable that facilitates, in each step of a given temporal discretization scheme, an effective splitting between the PDE part and the early exercise constraint.
This IT splitting technique has next been employed in~\cite{IT09} for the Heston PDCP and in~\cite{T08} for the Kou PIDCP.
For treating the integral part in the latter case, an iterative method has been applied that is similar to a fixed-point iteration.

Haentjens et al.~\cite{HH15,HHV10} considered the Heston PDCP and combined alternating direction implicit (ADI) schemes for directional splitting of PDEs with the IT splitting technique for the early exercise constraint, defining the so-called ADI-IT methods.
These methods were next studied in~\cite{HV17} for the one- and two-dimensional Black--Scholes PDCPs, where also a useful interpretation of this combined splitting approach was provided.
 
Complementary to this, the adaptation of ADI schemes to two-dimensional partial integro-differential equations (PIDEs) has recently been investigated by in 't Hout \& Toivanen~\cite{HT18}.
Here three novel adaptations of the well-established modified Craig--Sneyd (MCS) scheme~\cite{HW09} were analyzed and applied for the valuation of European options under the Bates model, where the mixed derivative term and the integral part are conveniently treated in an explicit fashion.

Boen \& in 't Hout~\cite{BH19} subsequently studied seven operator splitting schemes of both the implicit-explicit (IMEX) and the ADI kind in the application to the two-dimensional Merton PIDE, where the two-dimensional integral part is always handled explicitly. 
It was concluded that, among the schemes considered, the adaptation introduced in~\cite{HT18} of the MCS scheme where the integral part is dealt with in a two-step Adams--Bashforth manner is preferable.

Finally, Heidarpour-Dehkordi \& Christara~\cite{HC18} recently proposed an adaptation of the MCS scheme for the two-dimensional Black--Scholes PDCP by means of the penalty method and applied it to value American spread options. 

The main aim of the present paper is to introduce and investigate adaptations of IMEX and ADI operator splitting schemes for the two-dimensional Merton PIDCP. 
We extend all second-order schemes from~\cite{BH19} by using the IT splitting technique and next study their performance in ample numerical experiments.
In addition, we consider two penalty type methods that were proposed in \cite{HFL04} and \cite{HC18}.

This paper is organized as follows.
In Section~\ref{sec:PIDC} the PIDCP is given that holds for American-style option values under the two-asset Merton jump-diffusion model.
Section~\ref{spatdisc} describes the spatial discretization of this PIDCP. 
Section~\ref{tempdisc} deals with the temporal discretization of the obtained semidiscrete PIDCP and defines the adaptation of the pertinent IMEX and ADI schemes considered in~\cite{BH19} by means of IT splitting, which will be applied in an iterative manner.
Also the two penalty type methods from \cite{HFL04} and \cite{HC18} are formulated in this section.
In Section~\ref{num res} extensive numerical experiments are presented that yield detailed insight into the temporal convergence behaviour of all operator splitting methods from Section~\ref{tempdisc} and their relative performance. 
Section~\ref{sec:conc} gives our conclusions.

\clearpage
\setcounter{equation}{0}
\section{PIDC problem}\label{sec:PIDC}
The PIDE for the value $v = v(s_1,s_2,t)$ of a European-style option with maturity time $T>0$ under the two-asset Merton jump-diffusion model is 
given by
\begin{equation}\label{PIDE2D}
\dud{t} =  \mathcal{D} v + \mathcal{J} v,
\end{equation}
with differential and integral operators
\begin{align*}
\mathcal{D} v =&~\tfrac{1}{2} \sigma_1^2s_1^2\dudd{s_1} + \rho \sigma_1\sigma_2s_1s_2\duddm{s_1}{s_2} + \tfrac{1}{2} \sigma_2^2 s_2^2\dudd{s_2} + 
(r-\lambda \zeta_1) s_1\dud{s_1} + (r-\lambda\zeta_2) s_2 \dud{s_2} -(r+\lambda)v, \\
\mathcal{J} v =&~\lambda\int_0^{\infty}\int_0^{\infty} v(s_1y_1, s_2y_2,t) f(y_1,y_2) \dd y_1 \dd y_2.
\end{align*}
Here $0<t\leq T$ and $s_q>0$ ($q=1,2$) represents the price of asset $q$ at time $\tau = T-t$. Next,
\begin{itemize}
\item $r$ is the risk-free interest rate,
\item $\sigma_q$ ($q=1,2$) is the volatility of asset $q$, conditional on the event that no jumps occur,
\item $\rho$ is the correlation between the two underlying standard Brownian motions,
\item $\lambda$ is the jump intensity of the underlying Poisson arrival process, 
\item $\zeta_q$ ($q=1,2$) is the expected relative jump size for asset $q$.
\end{itemize}
Function $f$ is the probability density function of a bivariate lognormal distribution,

\begin{align*}
f(y_1,y_2) = \frac{1}{2\pi \delta_1\delta_2\sqrt{1-\rhoh^2}\,y_1y_2}\exp\left(-\frac{\left(\frac{\ln(y_1)-\gamma_1}{\delta_1}\right)^2+
\left(\frac{\ln(y_2)-\gamma_2}{\delta_2}\right)^2-2\rhoh\left(\frac{\ln(y_1)-\gamma_1}{\delta_1}\right)
\left(\frac{\ln(y_2)-\gamma_2}{\delta_2}\right)}{2(1-\rhoh^2)}\right)
\end{align*}
for $y_1>0$, $y_2>0$.
Here $\gamma_q$ and $\delta_q$ ($q=1,2$) and $\rhoh$ are given real constants that can be interpreted as, respectively, the mean and standard 
deviation and correlation of a bivariate normal distribution. It holds that 
\[
\zeta_q = e^{\gamma_q+\tfrac{1}{2} \delta_q^2}-1\quad (q=1,2).
\]
Let  $\phi$ denote the payoff function of the option. Then we have the initial condition
\begin{equation}\label{IC}
v(s_1,s_2,0) = \phi(s_1,s_2)
\end{equation}
for $s_1>0$, $s_2>0$.

With the above notations and initial condition, the value of an American-style option satisfies the following PIDCP,
\begin{align}\label{PIDC}
\begin{cases}
\displaystyle\dud{t}(s_1,s_2,t) \geq  \mathcal{D} v(s_1,s_2,t) + \mathcal{J} v(s_1,s_2,t),\\\\
v(s_1,s_2,t) \geq \phi(s_1,s_2),\\\\
\left(v(s_1,s_2,t)-\phi(s_1,s_2)\right)\left(\displaystyle\dud{t}(s_1,s_2,t)-\mathcal{D} v(s_1,s_2,t) - \mathcal{J} v(s_1,s_2,t)\right) = 0
\end{cases}
\end{align}
whenever $s_1>0$, $s_2>0$, $0<t\leq T$.
Boundary conditions are given by imposing \eqref{PIDC} for $s_1=0$ and $s_2=0$, respectively.

\setcounter{equation}{0}
\section{Spatial discretization}\label{spatdisc}
In this section we describe the spatial discretization of the PIDCP \eqref{PIDC}.
As the first step, the unbounded spatial domain is truncated to $[0,S_{\rm max}]\times[0,S_{\rm max}]$ with sufficiently large value 
$S_{\rm max}$.
On the two far sides $s_{1} = S_{\rm max}$ and $s_{2} = S_{\rm max}$ linear boundary conditions are taken, which are common in finance,
\begin{equation}\label{LBC}
\dudd{s_1} = 0~~(\textrm{if}~s_{1} = S_{\rm max}) \quad \mbox{ and } \quad  \dudd{s_2} = 0 ~~(\textrm{if}~s_{2} = S_{\rm max}).
\end{equation}

Let $K>0$ denote the strike price of the option.
As in \cite{HV17}, a smooth nonuniform Cartesian spatial grid is constructed such that $K$ falls midway two 
successive grid points in each direction. This turns out to be beneficial for accuracy.
Let parameter $d>0$ and fix a subinterval $[S_{\mathrm{left}},S_{\mathrm{right}}]$ of $[0,S_{\rm max}]$ such that 
$S_{\mathrm{left}}+S_{\mathrm{right}} = 2K < S_{\rm max}$. 
Define
\begin{align*}
\xi_{\min} &=  \sinh^{-1} \left(\frac{-S_{\mathrm{left}}}{d}\right),\\
\xi_{\mathrm{int}} &= \frac{S_{\mathrm{right}}-S_{\mathrm{left}}}{d},\\
\xi_{\max} &= \xi_{\mathrm{int}} + \sinh^{-1}\left(\frac{S_{\max}-S_{\mathrm{right}}}{d}\right).
\end{align*}
Let $q\in \{1,2\}$.
For any given integer $\nu = \nu_q \geq1$, let
\[
\Delta \xi = \Delta \xi_q = \frac{\xi_{\mathrm{int}}-2\xi_{\min}}{\nu},
\]
let $m = m_q >\nu$ be the smallest integer such that $m\Delta \xi \geq \xi_{\max} - \xi_{\min}$, reset $\xi_{\max}$ to 
$\xi_{\min} + m\Delta \xi$ and define equidistant points 
\[
\xi_j = \xi_{q,j} = \xi_{\min} + j \Delta \xi \quad (j=0,1,\ldots,m).
\] 
It is easily verified that $\xi_{\mathrm{int}}/2$ forms the middle of $[\xi_0, \xi_\nu]$.
Hence, it lies exactly halfway two successive $\xi$-grid points whenever $\nu$ is odd. 
The grid in the $q$-th spatial direction 
\[
0 = s_{q,0} < s_{q,1} < \ldots < s_{q, m}
\]
is then defined through the transformation 
\[
s_{q,j} = \psi(\xi_j) \quad (j=0,1,\ldots,m),
\]
where 
\[
\psi(\xi) = \begin{cases}
S_{\mathrm{left}} + d \sinh(\xi) & \textrm{for}\ \xi_{\min} \leq \xi \leq 0,\\
S_{\mathrm{left}} + d \xi & \textrm{for}\ 0 < \xi \leq \xi_{\mathrm{int}},\\
S_{\mathrm{right}} + d \sinh(\xi - \xi_{\mathrm{int}}) &\textrm{for}\ \xi_{\mathrm{int}} < \xi \leq \xi_{\max}.
\end{cases}
\]
This grid is uniform inside $[S_{\mathrm{left}},S_{\mathrm{right}}]$ and nonuniform outside, where the mesh width 
inside equals $d \Delta \xi$ and is always smaller than the mesh widths outside.
By construction, $K = \psi(\xi_\mathrm{int}/2)$ falls halfway two successive grid points in each spatial direction 
whenever $\nu$ is odd.  
In this paper we heuristically select $d=K/3$ and $S_{\mathrm{left}} = 0.8K$, $S_{\mathrm{right}}=1.2K$.
Further, $S_{\rm max}$ is reset to $\psi(\xi_{\max})$, which is slightly larger than the original value.

The discretization of the convection-diffusion-reaction part $\mathcal{D} v$ in \eqref{PIDC} is performed using 
finite differences. 
Let $u:[0,S_{\rm max}]\to\mathbb{R}$ be any given smooth function, let $0=s_0<s_1<\ldots<s_m=S_{\rm max}$ be any 
given smooth nonuniform unidirectional grid and define mesh widths $h_j = s_j - s_{j-1}$ ($1\le j\le m$). 
We consider the following second-order central finite difference formulas for convection and diffusion:
\begin{align*}
u'(s_j) & \approx \alpha_{j,-1}u(s_{j-1}) + \alpha_{j,0}u(s_j) + \alpha_{j,1}u(s_{j+1}),\\
\nonumber\\
u''(s_j) & \approx \beta_{j,-1}u(s_{j-1}) + \beta_{j,0} u(s_{j}) + \beta_{j,1}u(s_{j+1}),
\end{align*}
with
\begin{align*}
   \alpha_{j,-1} &=  \frac{-h_{j+1}}{h_j(h_j+h_{j+1})}  
 & \alpha_{j,0} &= \frac{h_{j+1}-h_j}{h_j h_{j+1}}  
 & \alpha_{j,1} &= \frac{h_j}{h_{j+1}(h_j+h_{j+1})},\\
   \beta_{j,-1} &= \frac{2}{h_j(h_j+h_{j+1})}   
 & \beta_{j,0} &= \frac{-2}{h_j h_{j+1}} 
 & \beta_{j,1} &= \frac{2}{h_{j+1}(h_j+h_{j+1})}
\end{align*}
for $1\le j\le m-1$.
In view of the degeneracy of $\mathcal{D} v$ at the zero boundaries, no discretization is required if $j=0$. 
If $j=m$, then we use the first-order backward finite difference formula for the first derivative and, by the linear
boundary condition \eqref{LBC}, the second derivative vanishes. 
Finally, the mixed derivative term in $\mathcal{D} v$ is handled by applying the finite difference formulas for the 
first derivative consecutively in the two spatial directions.

Let the vector $V(t) = (V_{0,0}(t), V_{1,0}(t),\ldots,V_{m_1-1,m_2}(t),V_{m_1,m_2}(t))^{\top}$ where entry $V_{i,j}(t)$ 
denotes the semidiscrete approximation to $v(s_{1,i},s_{2,j},t)$ for $0\leq i \leq m_1$, $0\leq j \leq m_2$.
The semidiscrete version of the convection-diffusion-reaction part $\mathcal{D} v$ can then be written as 
\[
A^{(D)} V(t)
\] 
with matrix
\begin{equation*}
A^{(D)} = A^{(M)} + A_1 + A_2,
\end{equation*}
where the matrix $A^{(M)}$ corresponds to the mixed derivative term and the matrix $A_q$ corresponds to all derivative 
terms in the $q$-th spatial direction ($q=1,2$). 
Further, the reaction term has been distributed equally across $A_1$ and $A_2$.

For the discretization of the integral part $\mathcal{J}v$ we consider a transformation to the log-price variable
$x_q = \ln(s_q)$ ($q=1,2$).
This yields a two-dimensional cross-correlation,
\begin{equation}\label{2Dcrocor}
\left(\overline{\mathcal{J}}\,\overline{v}\right)(x_1,x_2,t) := 
\lambda\int_{-\infty}^{\infty}\int_{-\infty}^{\infty} \overline{v}(z_1,z_2,t) \overline{f}(z_1-x_1,z_2-x_2)\dd z_1 \dd z_2 =
(\mathcal{J}v)(e^{x_1},e^{x_2},t).
\end{equation}
Here $\overline{v}(z_1,z_2,t)=v(e^{z_1},e^{z_2},t)$ and $\overline{f}(\eta_1,\eta_2) = 
f(e^{\eta_1},e^{\eta_2})e^{\eta_1}e^{\eta_2}$ is the probability density function of a bivariate normal distribution. 
Let $M_q$ be a given power of 2 such that $\Delta x_q = \ln(S_{\rm max})/M_q$ is smaller than the smallest mesh width 
in the nonuniform $\ln(s_q)$-grid ($q=1,2$).
Then the double integral \eqref{2Dcrocor} is approximated on the uniform Cartesian grid $(x_{1,k},x_{2,l}) = (k\Delta x_1,l\Delta x_2)$ 
for ${k=-M_1+1,\ldots,M_1, \ l=-M_2+1,\ldots,M_2}$ by
\[
\overline{J}_{k,l}(t) =\lambda \sum_{j=-M_2+1}^{M_2}\sum_{i=-M_1+1}^{M_1} \overline{V}_{i,j}(t)\overline{f}_{i-k,j-l} \Delta x_1 \Delta x_2.
\]
Here $\overline{V}_{i,j}(t) \approx \overline{v}(x_{1,i}, x_{2,j},t)$ and $\overline{f}_{i-k,j-l} = \overline{f}(x_{1,i}-x_{1,k},x_{2,j}-x_{2,l})$. 
Define the vectors
\begin{align*}
\overline{J}(t) =&\ (\overline{J}_{-M_1+1, -M_2+1}(t), \overline{J}_{-M_1+2, -M_2+1}(t), \ldots, \overline{J}_{M_1-1, M_2}(t), \overline{J}_{M_1, M_2}(t))^\top,\\
\overline{V}(t) =&\ (\overline{V}_{-M_1+1, -M_2+1}(t), \overline{V}_{-M_1+2, -M_2+1}(t), \ldots, \overline{V}_{M_1-1, M_2}(t), \overline{V}_{M_1, M_2}(t))^\top,
\end{align*}
then
\[
\overline{J}(t) =\overline{A}^{(J)} \ \overline{V}(t),
\]
where $\overline{A}^{(J)}$ is a given asymmetric two-level block-Toeplitz matrix of size $(4M_1M_2) \times (4M_1M_2)$. 

Matrix-vector products with asymmetric multilevel block-Toeplitz matrices can be computed efficiently using the FFT algorithm 
introduced by Barrowes, Teixeira \& Kong~\cite{BTK01}. 
This algorithm employs two FFTs and one inverse FFT, all of them one-dimensional, reducing the computational cost of the 
matrix-vector product to $\mathcal{O}(M_1M_2\log(M_1M_2))$. 
The algorithm from~\cite{BTK01} essentially embeds the block-Toeplitz matrix in a circulant matrix before applying the FFT, thus 
avoiding any wrap-around effect, which would have occurred when applying the FFT directly to the block-Toeplitz matrix.

Notice that the price grid on which the PIDCP \eqref{PIDC} is discretized is nonuniform and in general does not define 
a uniform log-price grid. 
We therefore bilinearly interpolate the option value approximations between the $s$- and $x$-grids immediately before 
and after application of the FFT algorithm of \cite{BTK01}, leading to the discretized version of $\mathcal{J}v$:
\[
A^{(J)} V(t)
\]
with given, fixed matrix $A^{(J)}$.
More details are given in~\cite{BH19}.

With the above discretizations of $\mathcal{D} v$ and $\mathcal{J}v$, the following semidiscrete PIDCP is obtained:
\begin{align}\label{PIDCdisc}
V'(t) \geq AV(t), \quad V(t) \geq V^0, \quad  (V(t) - V^0)^\top(V'(t)-AV(t))=0
\end{align}
for $0<t\leq T$. Here inequalities for vectors are to be understood componentwise. The matrix 
\[
A = A^{(D)} + A^{(J)} = A^{(M)} + A_1 + A_2 + A^{(J)}
\] 
and the vector $V^0$ is directly given by the payoff function $\phi$.

It can be readily seen that the maximum norm of the matrix $A^{(J)}$ is bounded by a moderate constant, independent of the spatial 
grid, whenever the jump intensity $\lambda$ is moderate. 
Hence, the dense matrix $A^{(J)}$ forms a nonstiff part of the semidiscrete system, whereas the sparse matrix $A^{(D)}$ constitutes a 
stiff part.

\setcounter{equation}{0}
\section{Temporal discretization}\label{tempdisc}
In the temporal discretization of the semidiscrete PIDCP \eqref{PIDCdisc} we employ two main techniques that have 
been proposed in the computational finance literature: the Ikonen--Toivanen (IT) splitting technique and the penalty method.
The IT splitting technique was first considered for American option valuation in \cite{IT04,IT09,T08} and the penalty method in 
\cite{FV02,ZFV98}.

\subsection{IT splitting} \label{IT schemes}
IT splitting is combined with a temporal discretization scheme for semidiscrete PIDEs, so as to yield a temporal discretization of 
semidiscrete PIDCPs. 
In \cite{BH19} a variety of contemporary operator splitting schemes has been studied for the temporal discretization of the semidiscrete 
two-dimensional Merton PIDE for European option values. 
Here, the integral part, which corresponds to the dense matrix $A^{(J)}$, is always conveniently treated in an explicit fashion.

In the following we define six temporal discretization methods for the semidiscrete two-dimensional Merton PIDCP \eqref{PIDCdisc} 
by combining six operator splitting schemes from \cite{BH19} with the IT splitting technique.
The pertinent operator splitting schemes are of the implicit-explicit (IMEX) and the alternating direction implicit (ADI) kind.
The combination of IT splitting with IMEX schemes was introduced in \cite{KL11} and next investigated in e.g.~\cite{HT16,ST12,STS14} for one- and two-dimensional PIDCPs.
The combination of this technique with ADI schemes was introduced in \cite{HHV10} and subsequently studied in \cite{HH15,HV17} for two-dimensional PDCPs (without integral term).

A novel feature of our present application of the IT splitting technique is that it is used in an iterative manner in each time step. 
This turns out to yield a significant enhancement, as will be shown in Section~\ref{num res} below.
We denote the number of iterations by $\kappa\ge 1$ and the iterated version of the IT splitting technique by IT$(\kappa)$.
For $\kappa=1$, this reduces to the original IT splitting.

Let integer $N\ge 1$ be given, step size $\Delta t = T/N$ and $\wV^0 = V^0$.
Each of the following methods defines an approximation $\wV^n$ to $V(t^n)$ at the temporal grid point $t^n = n \Delta t$ successively for $n=1,2,\ldots,N$.

\begin{enumerate}
\item CNFI-IT$(\kappa)$ method:\\[0.2cm]
This method can be viewed as a combination of the Crank--Nicolson scheme with fixed-point iteration (CNFI) for PIDEs, proposed by Tavella \& Randall \cite{TR00book} and d'Halluin, Forsyth \& Vetzal \cite{HFV05}, with the IT$(\kappa)$ splitting technique:
\begin{align}\label{CNFI IT}
&\begin{cases}
    \left(I-\tfrac{1}{2}\Delta t A^{(D)}\right)Z_k = \left(I + \tfrac{1}{2}\Delta t A^{(D)}\right)\wV^{n-1} +  \tfrac{1}{2} \Delta t A^{(J)} (\widehat{Z}_{k-1}+\wV^{n-1}) + \Delta t \lambda_{k-1},\\[0.3cm]
    \widehat{Z}_k= \max\left\{Z_k-\Delta t \lambda_{k-1}, V^0\right\},\quad
    \lambda_k = \max\left\{0,\lambda_{k-1}+(V^0-Z_k)/\Delta t\right\},\\[0.3cm]
    \text{for} \ k=1,2,\ldots,\kappa \text{ and } \wV^n = \widehat{Z}_{\kappa},\quad \widehat{\lambda}^n = \lambda_{\kappa}.
\end{cases}
\end{align}
Here, and in the subsequent methods, $\widehat{Z}_0 = \wV^{n-1}$ and $\lambda_0=\widehat{\lambda}^{n-1}$ with $\widehat{\lambda}^0$ taken to be the zero vector. 
Clearly, in each time step, $\kappa$ matrix-vector products with the matrix $A^{(J)}$ occur. 
\end{enumerate}

\begin{enumerate}\setcounter{enumi}{1}   
\item IETR-IT$(\kappa)$ method:\\[0.2cm]
The IETR scheme has been considered in \cite{H17book} for PIDEs. It treats the convection-diffusion-reaction part using the implicit trapezoidal rule 
(Crank--Nicolson) and the integral part using the explicit trapezoidal rule. Combining this scheme with the IT($\kappa$) splitting technique yields:
\begin{align}\label{IETR IT}
    &\begin{cases}
    Y_0 = \wV^{n-1} + \Delta t\, (A^{(D)}+A^{(J)})\wV^{n-1} + \Delta t \lambda_{k-1},\\
    \bar{Y}_0 = Y_0+\tfrac{1}{2} \Delta t A^{(J)} \big(Y_0-\wV^{n-1}\big),\\
    Y_1 = \bar{Y}_0+\tfrac{1}{2} \Delta t A^{(D)} \big(Y_1-\wV^{n-1}\big),\\[0.3cm]
    Z_k  = Y_1,\\
    \lambda_k = \max\left\{0,\lambda_{k-1}+(V^0-Z_k)/\Delta t\right\}, \\[0.3cm]
    \text{for} \ k=1,2,\ldots,\kappa \text{ and } \\[0.3cm] 
    \wV^n = \max\left\{Z_\kappa-\Delta t \lambda_{\kappa-1}, V^0\right\},\\ 
    \widehat{\lambda}^n\, = \lambda_{\kappa}.
\end{cases}
\end{align}
Noticing that $A^{(J)}\wV^{n-1}$ can be computed upfront, the IETR-IT($\kappa$) method requires $\kappa+1$ matrix-vector products with $A^{(J)}$ per time step.

\item CNAB-IT$(\kappa)$ method:\\[0.2cm]
The CNAB scheme was proposed for PIDEs by Salmi \& Toivanen \cite{ST14}.
It again treats the convection-diffusion-reaction part by the Crank--Nicolson scheme, but the integral part is now handled in a 
two-step Adams--Bashforth manner.
In Salmi, Toivanen \& von Sydow \cite{STS14} the CNAB scheme has been combined with IT splitting and applied to two-dimensional PIDCPs.
The extension by IT($\kappa$) splitting yields:
\begin{align}\label{CNAB IT}
    &\begin{cases}
     \left(I-\tfrac{1}{2} \Delta t A^{(D)}\right)Z_k = \left(I+\tfrac{1}{2} \Delta t A^{(D)}\right)\wV^{n-1} + \tfrac{1}{2}\Delta t A^{(J)}(3\wV^{n-1}-\wV^{n-2}) + \Delta t \lambda_{k-1},\\[0.3cm]
    \lambda_k = \max\left\{0,\lambda_{k-1}+(V^0-Z_k)/\Delta t\right\}, \\[0.3cm]
    \text{for} \ k=1,2,\ldots,\kappa \text{ and } \\[0.3cm] 
    \wV^n = \max\left\{Z_\kappa-\Delta t \lambda_{\kappa-1}, V^0\right\},\\ 
    \widehat{\lambda}^n\, = \lambda_{\kappa}.
\end{cases}
\end{align}
Clearly, each time step of this method requires just one matrix-vector product with $A^{(J)}$.
\end{enumerate}

The temporal discretization schemes that underly the three methods \eqref{CNFI IT}, \eqref{IETR IT}, \eqref{CNAB IT} are all of the IMEX kind, 
making use of the operator splitting $A = A^{(D)} + A^{(J)}$.
The next three methods are based upon ADI schemes, which employ the additional splitting $A^{(D)} = A^{(M)} + A_1 + A_2$.
\begin{enumerate}\setcounter{enumi}{3}   
\item MCS-IT$(\kappa)$ method:\\[0.2cm]
The modified Craig--Sneyd (MCS) scheme was introduced for PDEs containing mixed derivative terms by in 't Hout \& Welfert \cite{HW09}.
Its direct adaptation to PIDEs has been investigated by in 't Hout \& Toivanen \cite{HT18}.
Complementary to this, the MCS scheme has been combined with IT splitting by Haentjens et al. \cite{HH15,HHV10}, defining 
the MCS-IT method for PDCPs (without integral part).
The following MCS-IT$(\kappa)$ method generalizes all of these:
\begin{align}\label{MCS IT}
    &\begin{cases}
    Y_0 = \wV^{n-1} + \Delta t\, (A^{(D)}+A^{(J)})\wV^{n-1} + \Delta t \lambda_{k-1},\\
    Y_j = Y_{j-1} + \theta\Delta t A_j(Y_j-\wV^{n-1}) \quad (j=1,2),\\
    \bar{Y}_0 = Y_0 + \theta \Delta t\, (A^{(M)}+A^{(J)})(Y_2 - \wV^{n-1}),\\
    \widetilde{Y}_0 = \bar{Y}_0 + (\tfrac{1}{2}-\theta)\Delta t\, (A^{(D)}+A^{(J)})(Y_2-\wV^{n-1}),\\
    \widetilde{Y}_j = \widetilde{Y}_{j-1} + \theta \Delta t A_j(\widetilde{Y}_j - \wV^{n-1}) \quad (j=1,2),\\[0.3cm]
    Z_k = \widetilde{Y}_2,\\
    \lambda_k = \max\left\{0,\lambda_{k-1}+(V^0-Z_k)/\Delta t\right\}, \\[0.3cm]
    \text{for} \ k=1,2,\ldots,\kappa \text{ and } \\[0.3cm] 
    \wV^n = \max\left\{Z_\kappa-\Delta t \lambda_{\kappa-1}, V^0\right\},\\ 
    \widehat{\lambda}^n\, = \lambda_{\kappa}.
\end{cases}
\end{align}
We make the standard choice $\theta = \frac{1}{3}$, which is prompted by stability and accuracy results in the literature for two-dimensional problems (see e.g. \cite{HM11,HT18,HW09,HW16}). 
Since $A^{(J)}\widehat{V}^{n-1}$ can be computed in advance and the explicit stages $\bar{Y}_0$, $\widetilde{Y}_0$ can be merged, the integral part is evaluated $\kappa+1$ times per time step.

It is easily verified that in \eqref{MCS IT} the integral part and mixed derivative term are both handled by the explicit trapezoidal rule.
The implicit stages $Y_j$, $\widetilde{Y}_j$ (for $j=1,2$) are often called stabilizing corrections and are unidirectional.
The pertinent linear systems for these stages are tridiagonal and can therefore be solved very efficiently by using an a priori $LU$ factorization.\\

\clearpage
\item MCS2-IT$(\kappa)$ method:\\[0.2cm]
An alternative adaptation of the MCS scheme to PIDEs has been proposed in \cite{HT18} where the integral part is treated in a two-step Adams--Bashforth manner. Combining this adaptation with IT($\kappa$) splitting leads to:
\begin{align}\label{MCS2 IT}
    &\begin{cases}
    X_0 = \wV^{n-1}+\Delta t A^{(D)}\wV^{n-1} + \Delta t \lambda_{k-1},\\
    Y_0 = X_0 + \tfrac{1}{2}\Delta t A^{(J)}(3\wV^{n-1}-\wV^{n-2}),\\
    Y_j = Y_{j-1}+\theta\Delta t A_j(Y_j-\wV^{n-1})\quad (j=1,2),\\
    \bar{Y}_0 = Y_0 + \theta\Delta t A^{(M)}(Y_2-\wV^{n-1}),\\
    \widetilde{Y}_0 = \bar{Y}_0 + (\tfrac{1}{2}-\theta)\Delta t A^{(D)}(Y_2-\wV^{n-1}),\\
    \widetilde{Y}_j = \widetilde{Y}_{j-1} + \theta \Delta t A_j(\widetilde{Y}_j-\wV^{n-1})\quad (j=1,2),\\[0.3cm]
    Z_k = \widetilde{Y}_2,\\
    \lambda_k = \max\left\{0,\lambda_{k-1}+(V^0-Z_k)/\Delta t\right\}, \\[0.3cm]
    \text{for} \ k=1,2,\ldots,\kappa \text{ and } \\[0.3cm] 
    \wV^n = \max\left\{Z_\kappa-\Delta t \lambda_{\kappa-1}, V^0\right\},\\ 
    \widehat{\lambda}^n\, = \lambda_{\kappa}.
\end{cases}
\end{align}
We choose again $\theta = \frac{1}{3}$.
It is clear that the MCS2-IT$(\kappa)$ method requires only one evaluation of the integral part per time step.

\item SC2A-IT$(\kappa)$ method:\\[0.2cm]
The stabilizing correction two-step Adams-type (SC2A) scheme is a prominent member of the class of stabilizing correction multistep methods that has been investigated by Hundsdorfer \& in 't Hout \cite{HH18} for the numerical solution of PDEs.
Its direct adaptation to PIDEs has been studied by Boen \& in 't Hout \cite{BH19} and treats the integral part and mixed derivative term jointly in a two-step Adams--Bashforth fashion. 
The combination with IT($\kappa$) splitting yields:
\begin{align}\label{SC2A IT}
    &\begin{cases}
    X_0 = \wV^{n-1} + \Delta t\, (A_1+A_2)\sum_{i=1}^2b^0_i\wV^{n-i} + \Delta t \lambda_{k-1},\\
    Y_0 =  X_0 + \Delta t\, (A^{(M)}+A^{(J)})\sum_{i=1}^2b^1_i\wV^{n-i},\\
    Y_j = Y_{j-1} + \theta \Delta t A_j(Y_j-\wV^{n-1})\quad (j=1,2),\\[0.3cm]
    Z_k = Y_2,\\
    \lambda_k = \max\left\{0,\lambda_{k-1}+(V^0-Z_k)/\Delta t\right\}, \\[0.3cm]
    \text{for} \ k=1,2,\ldots,\kappa \text{ and } \\[0.3cm] 
    \wV^n = \max\left\{Z_\kappa-\Delta t \lambda_{\kappa-1}, V^0\right\},\\ 
    \widehat{\lambda}^n\, = \lambda_{\kappa}.
\end{cases}
\end{align}
Here the coefficients are $(b^1_1,b^1_2) = \left(\frac{3}{2},-\half\right)$ and $(b^0_1,b^0_2) = \left(\frac{3}{2}-\theta, - \half+\theta\right)$ and, following \cite{HH18}, we select $\theta = \frac{3}{4}$.
The SC2A-IT($\kappa$) method also requires just one evaluation of the integral part per time step.
\end{enumerate}
\vskip0.3cm
For each of the methods \eqref{CNFI IT}--\eqref{SC2A IT}, the underlying IMEX or ADI scheme has order of consistency equal to two for fixed, 
nonstiff systems of ordinary differential equations (ODEs), provided $\kappa \ge 2$ for method \eqref{CNFI IT}. 

In view of the nonsmoothness of the initial (payoff) function $\phi$, the first two time steps of each of the six methods above are replaced by 
four damping steps with step size $\Delta t/2$ using the backward Euler scheme with fixed-point iteration and IT$(\kappa)$ splitting, that is, 
the BEFI-IT$(\kappa)$ method.
With the full step size $\Delta t$, this method reads
\begin{align*}
    &\begin{cases}
    \left(I-\Delta t A^{(D)}\right)Z_k = \wV^{n-1} + \Delta t A^{(J)} \widehat{Z}_{k-1}+\Delta t \lambda_{k-1},\\[0.3cm]
    \widehat{Z}_k= \max\left\{Z_k-\Delta t \lambda_{k-1}, V^0\right\},\quad
    \lambda_k = \max\left\{0,\lambda_{k-1}+(V^0-Z_k)/\Delta t\right\},\\[0.3cm]
    \text{for} \ k=1,2,\ldots,\kappa \text{ and } \wV^n = \widehat{Z}_{\kappa},\quad \widehat{\lambda}^n = \lambda_{\kappa}.
\end{cases}
\end{align*}

At present the convergence theory for time stepping methods based on IT splitting is still under development.
In \cite{HH15} a useful relevant result was proved for the BE-IT method applied to PDCPs.
This result generalizes straightforwardly to the BEFI-IT$(1)$ method applied to PIDCPs,
\begin{align}\label{IMEX Euler IT}
&\begin{cases}
(I-\Delta t A^{(D)})Z = (I + \Delta t A^{(J)}) \widehat{V}^{n-1} +\Delta t \widehat{\lambda}^{n-1},\\[0.3cm]
\widehat{V}^n = \max\{Z-\Delta t \widehat{\lambda}^{n-1}, V^0\}, \quad \widehat{\lambda}^n = \max\{0, \widehat{\lambda}^{n-1} 
+ (V^0-Z)/\Delta t\}.
\end{cases}
\end{align}
The above can be regarded as the IMEX Euler-IT method.
The corresponding method for PIDCPs without IT splitting is the IMEX Euler method, which can be written as
\begin{align} \label{IMEX Euler}
&\begin{cases}
(I-\Delta t A^{(D)}) V^n = (I + \Delta t A^{(J)}) V^{n-1} +\Delta t \lambda^{n},\\[0.3cm]
V^n \geq V^0, \quad \lambda^{n}\geq 0, \quad (V^n-V^0)^{\top} \lambda^n = 0.
\end{cases}
\end{align}

Let $p=(m_1+1)(m_2+1)$.
For any given diagonal matrix $\mathcal{D}\in \R^{p\times p}$ with positive diagonal entries, define the scaled 
inner product by
\[
\langle U , V \rangle_\mathcal{D} = V^{\top} \mathcal{D}\, U
~~{\rm whenever}~U, V \in \R^p
\]
and let $\| \cdot \|_\mathcal{D}$ denote both the induced vector and matrix norms.
The following theorem forms a direct generalization of \cite[Thm.~3.1]{HH15} and is stated without proof.\\

\begin{theorem}\label{thm}
Consider the processes {\rm (\ref{IMEX Euler IT})} and {\rm (\ref{IMEX Euler})}.
Assume there exists a positive diagonal matrix $\mathcal{D}$ such that
\begin{equation*}\label{assA}
\mathcal{D}A^{(D)} + (A^{(D)})^{\top} \mathcal{D} ~{\it is~negative~semidefinite}.
\end{equation*}
Assume there are real constants $\mu$, $\nu$ independent of the spatial and temporal grids such that 
\begin{equation*}
\| A^{(J)} \|_\mathcal{D} \le \mu
\end{equation*}
and
\begin{equation*}\label{lambda}
\|\lambda^1\|_\mathcal{D} + \sum^{N}_{n=2} \|\lambda^n - \lambda^{n-1}\|_\mathcal{D} \leq \nu.
\end{equation*}
Then
\begin{equation*}
\max_{1\le n\le N} \|V^n - \widehat{V}^n\|_\mathcal{D} \leq \nu e^{\mu T} \, \Delta t
\end{equation*}
whenever $\Delta t = T/N$, integer $N\ge 1$.
\end{theorem}
\vskip0.2cm
\noindent

Theorem~\ref{thm} yields the useful result that the sequence $\{ \widehat{V}^n \}$ generated by (\ref{IMEX Euler IT}) is ${\cal O}(\Delta t)$ close 
to the sequence $\{ V^n \}$ defined by the basic method (\ref{IMEX Euler}).
For further results and a discussion of the assumptions in this theorem, we refer to \cite{B19book,HH15}.
   
\subsection{Penalty method}
We next consider the penalty method in combination with temporal discretization schemes for PIDEs, defining
temporal discretization methods for PIDCPs. 

\begin{enumerate}\setcounter{enumi}{6}   
\item CNFI-P method:\\[0.2cm]
This method forms a combination of the CNFI scheme with the penalty method and has been introduced by d'Halluin, Forsyth \& Labahn \cite{HFL04}:
\begin{align}\label{theta-P}
 \begin{cases}
 (I-\tfrac{1}{2} \Delta t^n A^{(D)} + P_{k-1}) Z_{k} = (I+\tfrac{1}{2}\Delta t^n A^{(D)})\wV^{n-1} + \tfrac{1}{2}\Delta t^n A^{(J)} (Z_{k-1} + \wV^{n-1})+ P_{k-1}V^0\\[0.3cm] 
 \text{for } k=1,2,\ldots,\kappa 
 \text{ and } \wV^n = Z_\kappa.
 \end{cases}
 \end{align}
Here $Z_0 = \wV^{n-1}$ is the starting value for the penalty iteration and $P_{k-1}$ denotes the diagonal matrix with $l$-th diagonal entry
\begin{align*}
\left(P_{k-1}\right)_{l,l} = \begin{cases}
\textit{Large} & \text{if }(Z_{k-1})_l < (V^0)_l,\\
0 & \text{otherwise}.
\end{cases}
\end{align*}
The common convergence criterion is
\begin{align}\label{convcrit}
\max_l \frac{|(Z_\kappa)_l-(Z_{\kappa-1})_l|}{\max\{1,|(Z_\kappa)_l|\}}< \textit{tol}.
\end{align}
We choose as natural values $\textit{tol} = 10^{-7}$ and $\textit{Large} = 10^7$. 

There are $\kappa$ matrix-vector products with $A^{(J)}$ in each time step of the CNFI-P method.
Notice that $\kappa$ depends on the time step number in view of the dynamic convergence criterion \eqref{convcrit}.

The CNFI-P method is applied with suitable nonuniform temporal grid points, which has been shown in e.g.~\cite{HFL04, FV02} to 
improve the temporal convergence behaviour.
Following \cite{IT09,RW14} we take
\begin{equation}\label{nonunitemp}
t^n = \left(\frac{n}{N}\right)^2T \quad (n=0,1,2,\ldots,N)
\end{equation}
and set $\Delta t^n = t^n - t^{n-1}$.

\item MCS-P method: \\[0.2cm]
Heidarpour-Dehkordi \& Christara \cite{HC18} combined the MCS scheme with the penalty method, defining the MCS-P method for PDCPs. 
We consider here its direct adaptation to PIDCPs:
\begin{align}\label{MCS P}
    \begin{cases}
    Y_0 = \wV^{n-1} + \Delta t\, (A^{(D)}+A^{(J)})\wV^{n-1},\\
    Y_j = Y_{j-1} + \theta\Delta t A_j(Y_j-\wV^{n-1}) \quad (j=1,2),\\
    \bar{Y}_0 = Y_0 + \theta \Delta t\, (A^{(M)}+A^{(J)})(Y_2 - \wV^{n-1}),\\
    \widetilde{Y}_0 = \bar{Y}_0 + (\tfrac{1}{2}-\theta)\Delta t\, (A^{(D)}+A^{(J)})(Y_2-\wV^{n-1}),\\
    \widetilde{Y}_1 = \widetilde{Y}_{0} + \theta \Delta t A_1(\widetilde{Y}_1 - \wV^{n-1}),\\[0.3cm]
    \text{and next} \\[0.3cm]
    (I-\theta \Delta t A_2 + P_{k-1} )Z_{k} = \widetilde{Y}_{1} - \theta \Delta t A_2 \wV^{n-1} + P_{k-1}V^0\\[0.3cm]
   \text{for} \ k=1,2,\ldots,\kappa \text{ and } \wV^n = Z_{\kappa}.
    \end{cases}
\end{align}
We select again $\theta = \frac{1}{3}$.
Clearly, the penalty iteration is introduced in the last implicit stage of the MCS scheme.
The same starting value $Z_0$, penalty matrix $P_{k-1}$, penalty factor $\textit{Large}$, convergence criterion \eqref{convcrit}
and tolerance $\textit{tol}$ are employed as for the CNFI-P method.
Since the penalty iteration does not involve the matrix $A^{(J)}$, the number of matrix-vector products with this matrix per time 
step of \eqref{MCS P} is equal to two.
\end{enumerate}

The first two time steps of the CNFI-P and MCS-P methods are replaced by four damping steps with half the original step size(s) 
applying the BEFI scheme combined with the penalty method, that is, the BEFI-P method. 

\setcounter{equation}{0}
\section{Numerical results}\label{num res}
\begin{figure}[h!]
    \centering
    \hspace{-0.5cm}\includegraphics[trim={1cm 0cm 1cm 0cm},clip,
    width =0.51\textwidth]{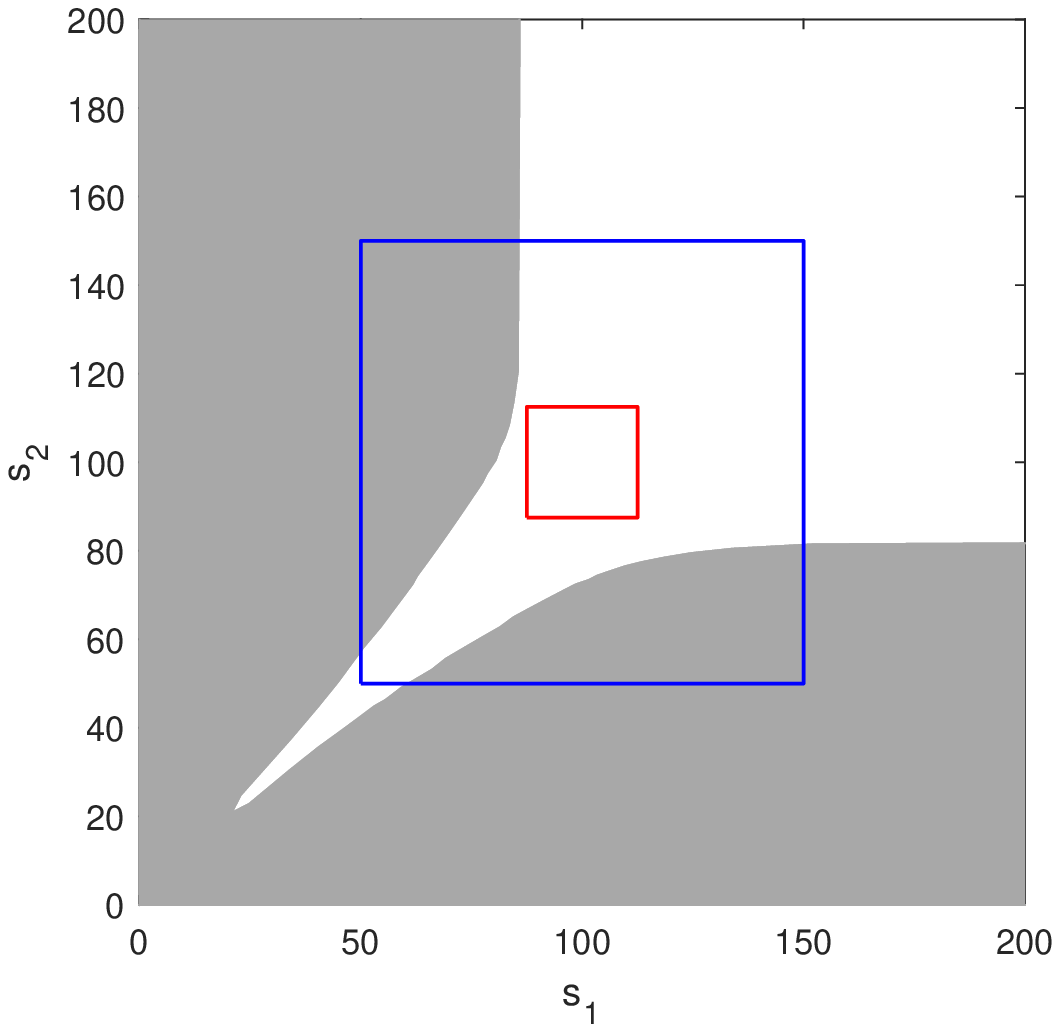}
    \includegraphics[trim={1cm 0cm 1cm 0cm},clip,
    width =0.51\textwidth]{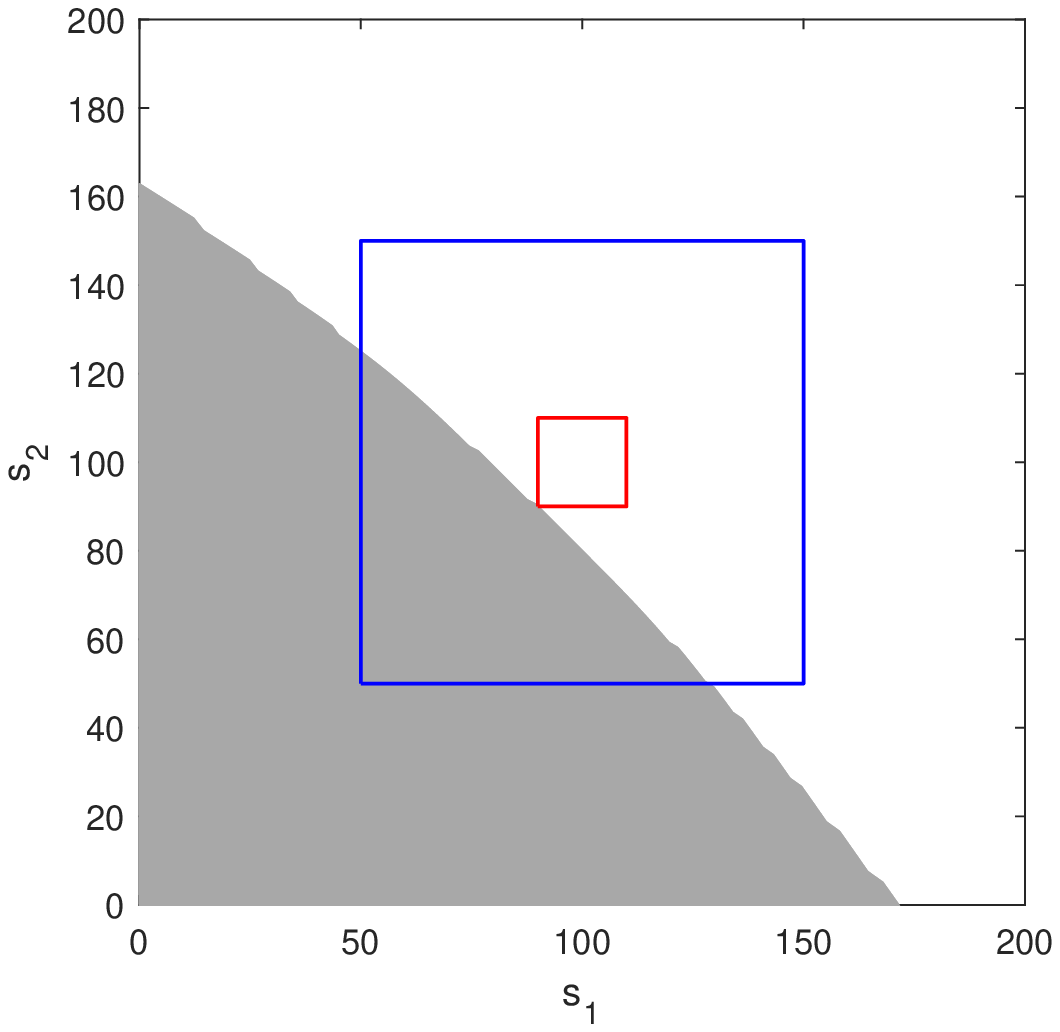}\\
    \hspace{-0.5cm}\includegraphics[trim={1cm 0cm 1cm 0cm},clip,
    width =0.51\textwidth]{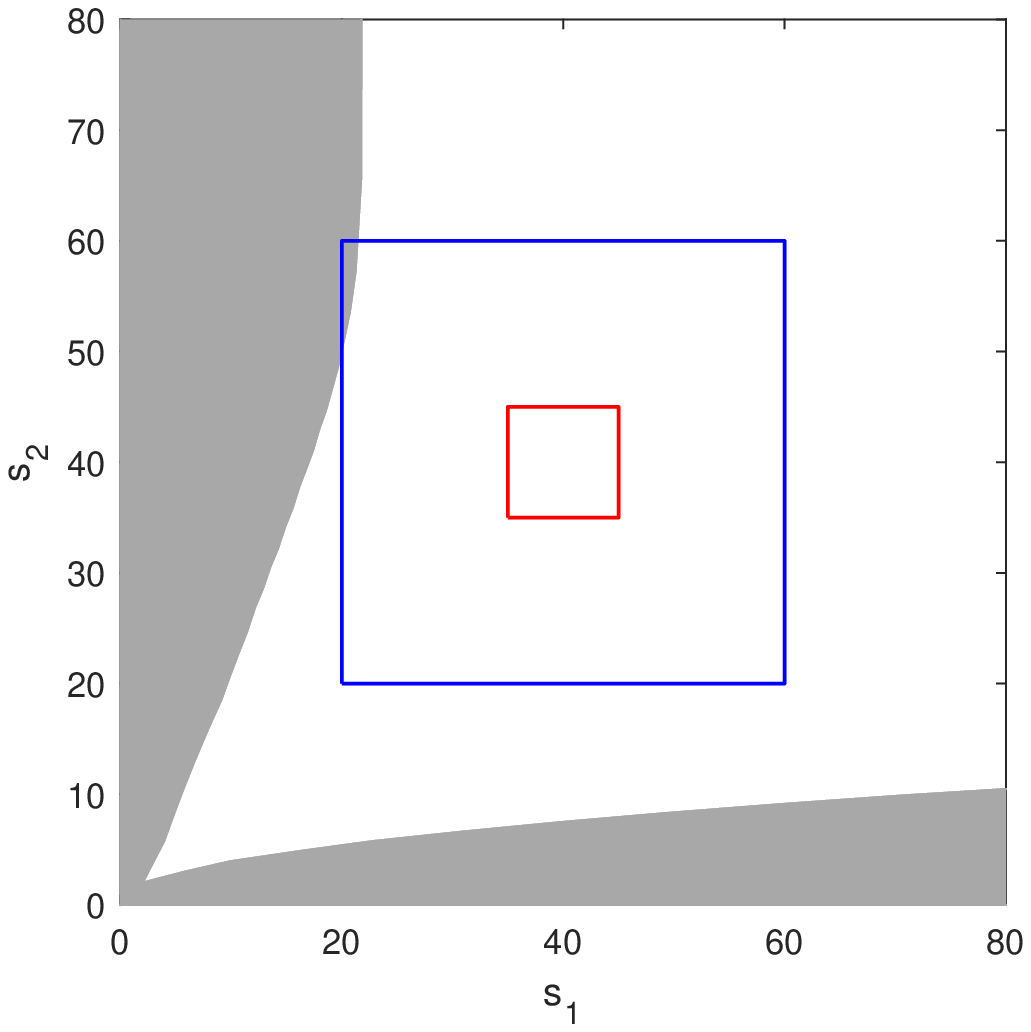}
    \includegraphics[trim={1cm 0cm 1cm 0cm},clip,
    width =0.51\textwidth]{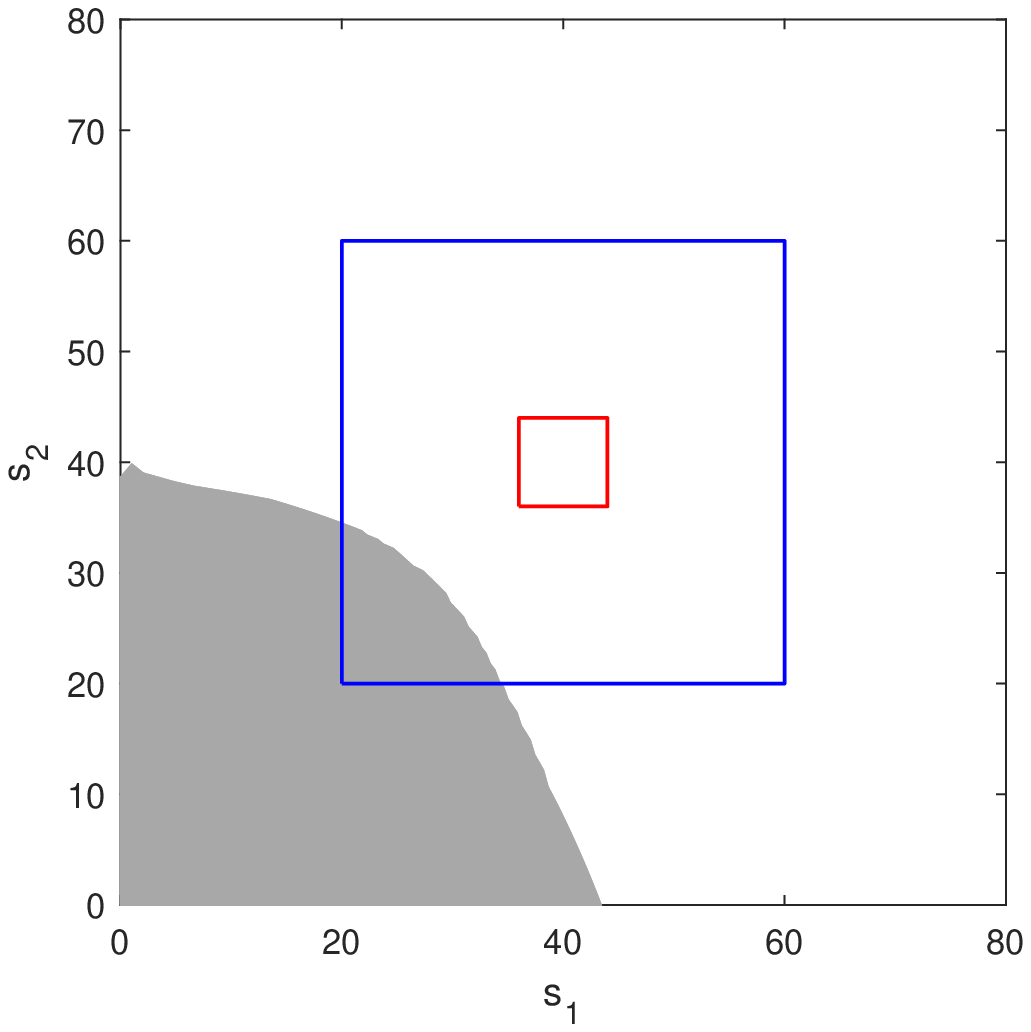}\\
    \hspace{-0.5cm}\includegraphics[trim={1cm 0cm 1cm 0cm},clip,
    width =0.51\textwidth]{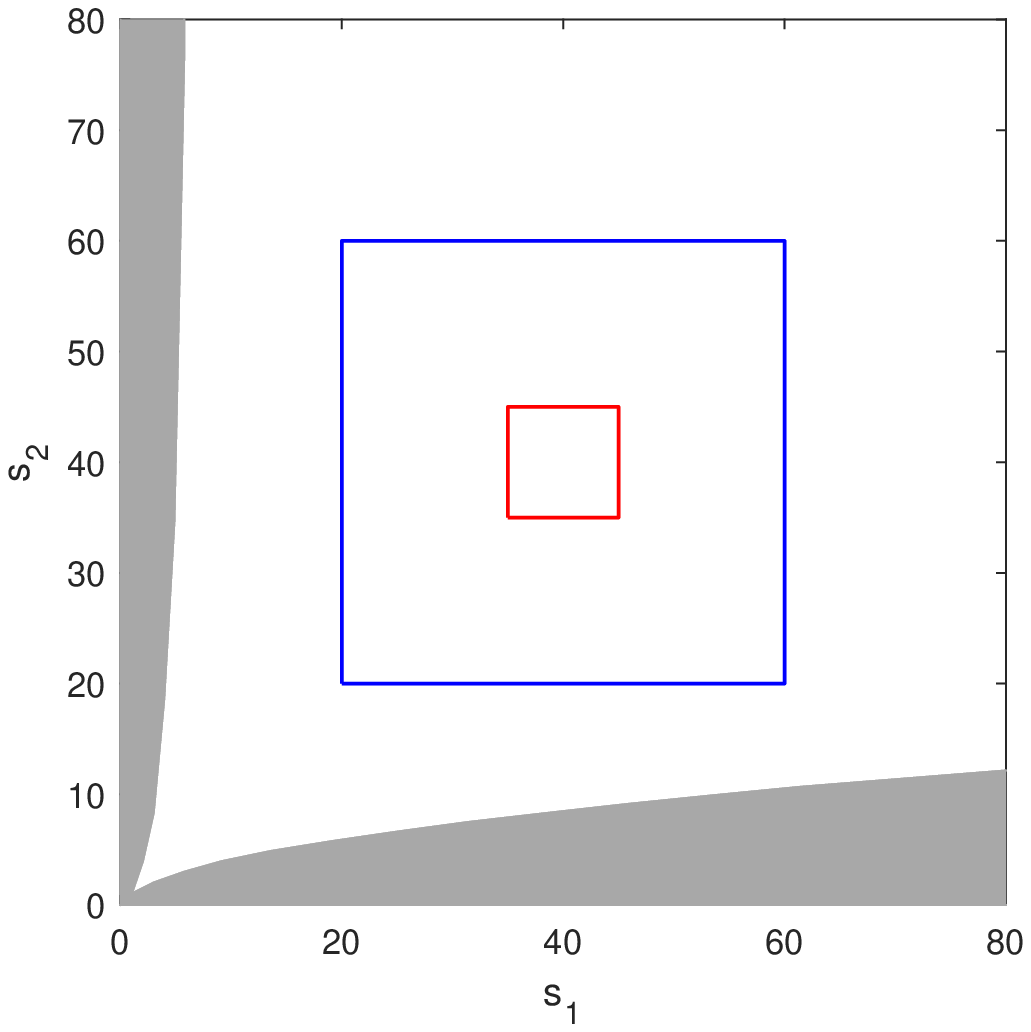}
    \includegraphics[trim={1cm 0cm 1cm 0cm},clip,
    width =0.51\textwidth]{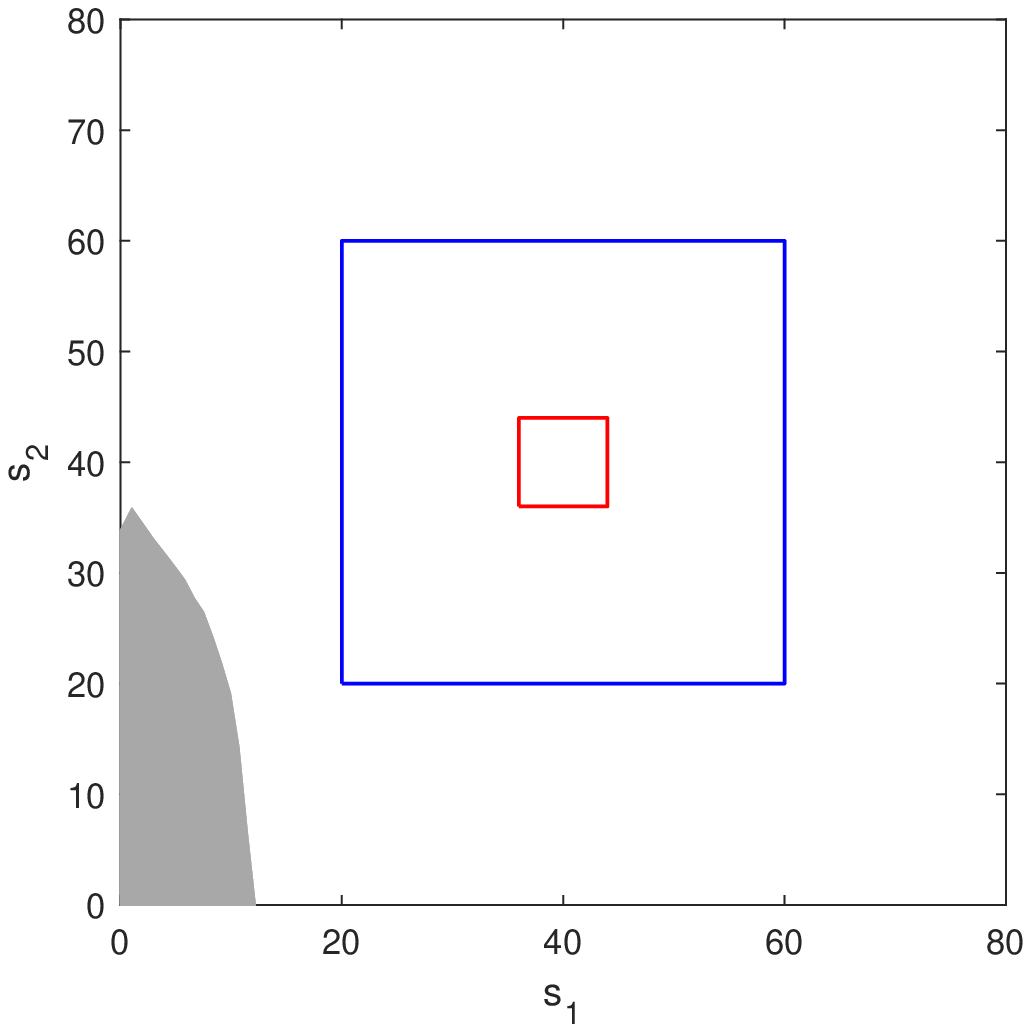}\\
    \caption{Early exercise regions for the American put-on-the-min (left) and put-on-the-average (right) options
    under the three parameter sets from 
    Table~\ref{tabpars} together with the regions of interest $[(1/2)K,(3/2)K]$ (blue) and $[(7/8)K,(9/8)K]$, respectively $[(9/10)K,(11/10)K]$ (red).}
    \label{EERs}
\end{figure}

\begin{figure}[h!]
    \centering
    \hspace{-0.5cm}\includegraphics[trim={0cm 0cm 0cm 0cm},clip,
    width =0.51\textwidth]{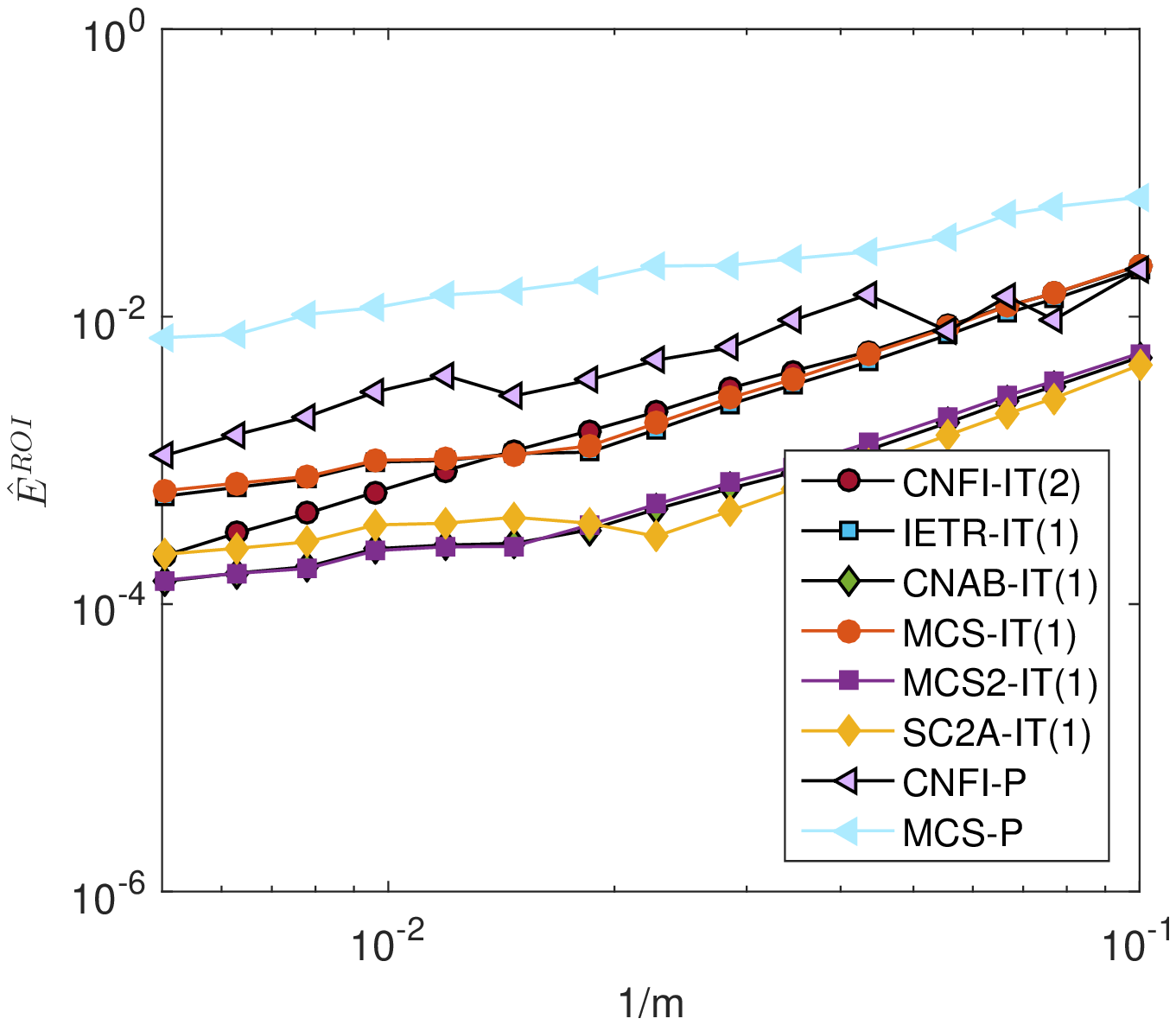}
    \includegraphics[trim={0cm 0cm 0cm 0cm},clip,
    width =0.51\textwidth]{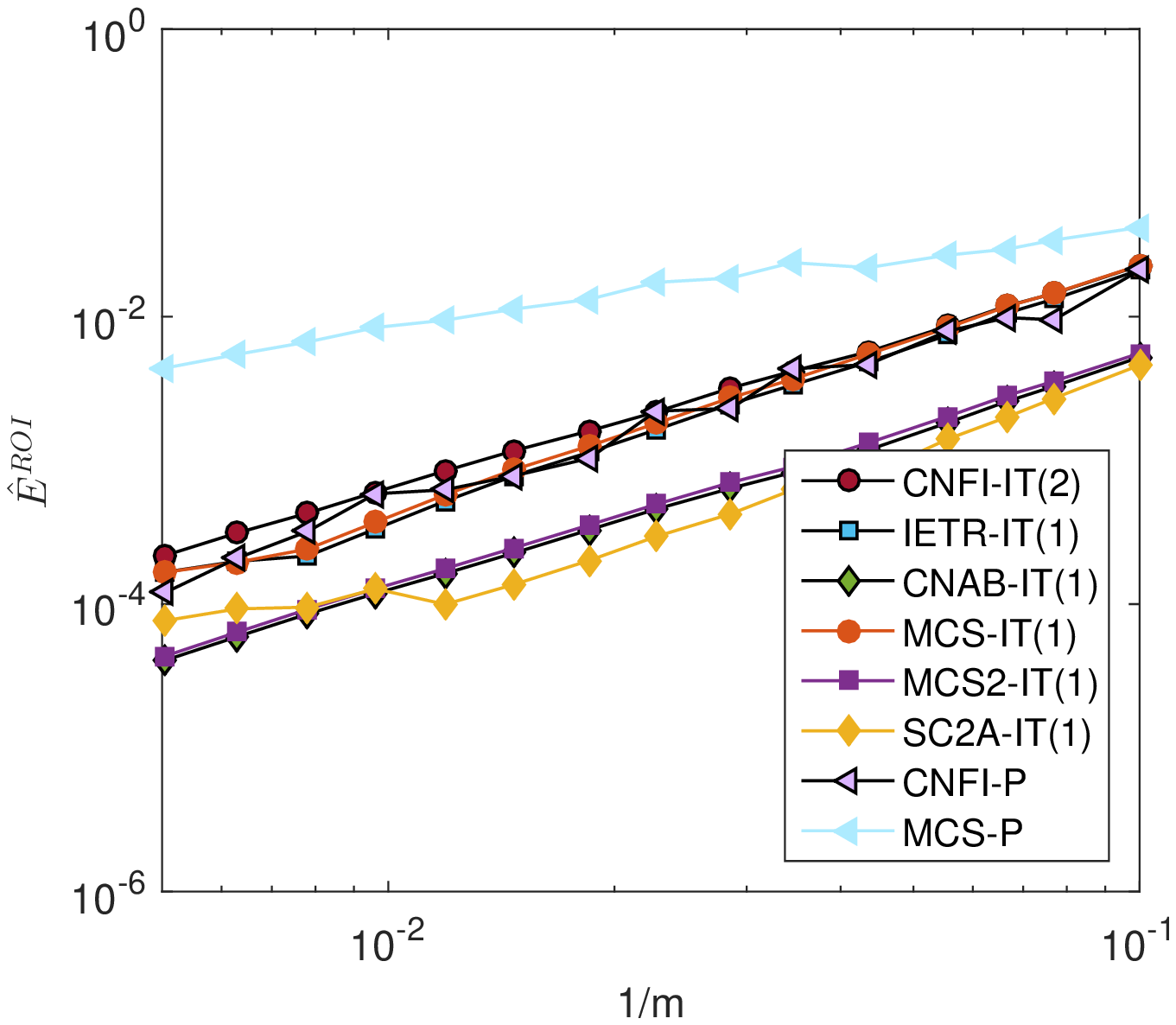}\\
    \hspace{-0.5cm}\includegraphics[trim={0cm 0cm 0cm 0cm},clip,
    width =0.51\textwidth]{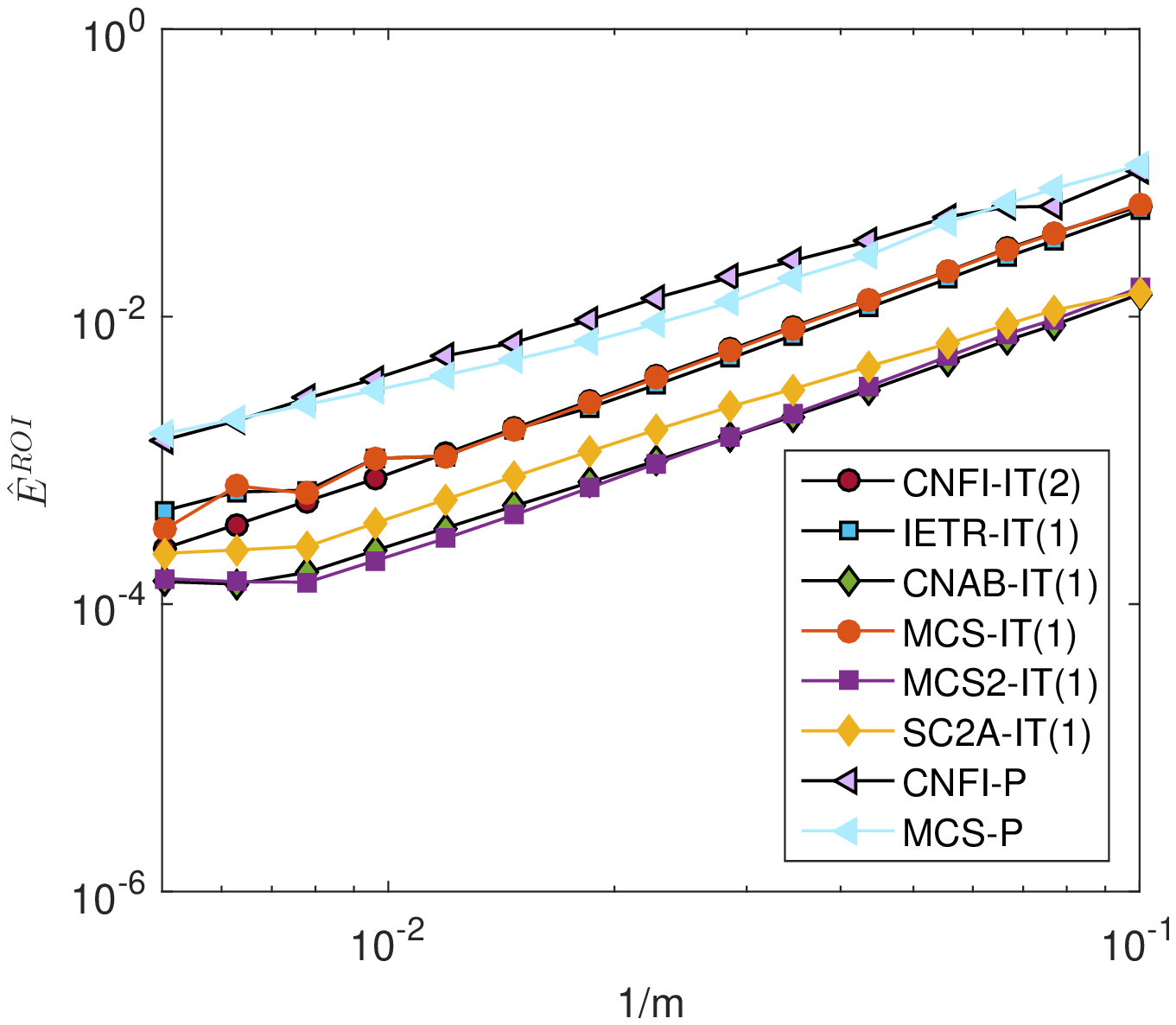}
    \includegraphics[trim={0cm 0cm 0cm 0cm},clip,
    width =0.51\textwidth]{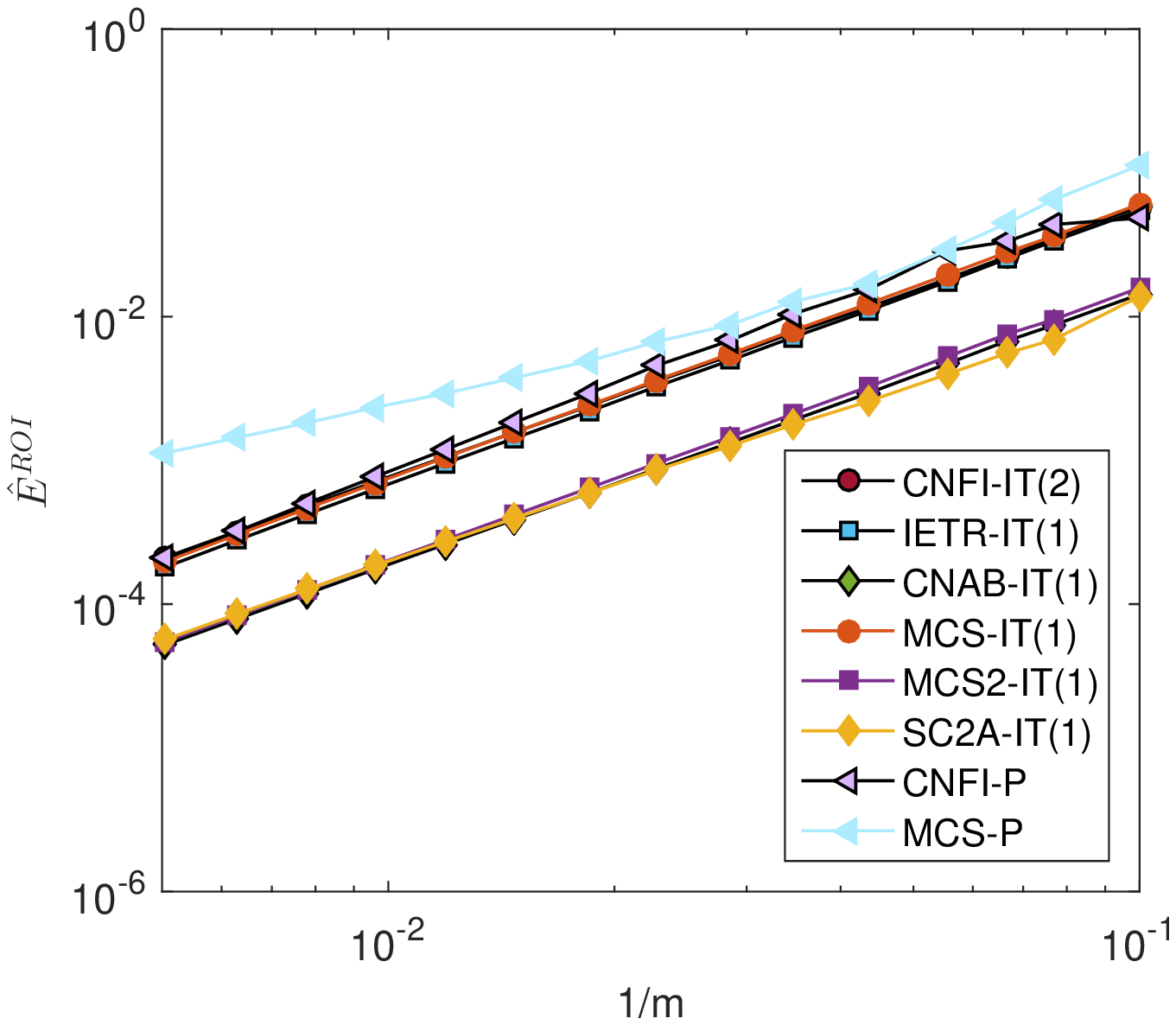}\\
    \hspace{-0.5cm}\includegraphics[trim={0cm 0cm 0cm 0cm},clip,
    width =0.51\textwidth]{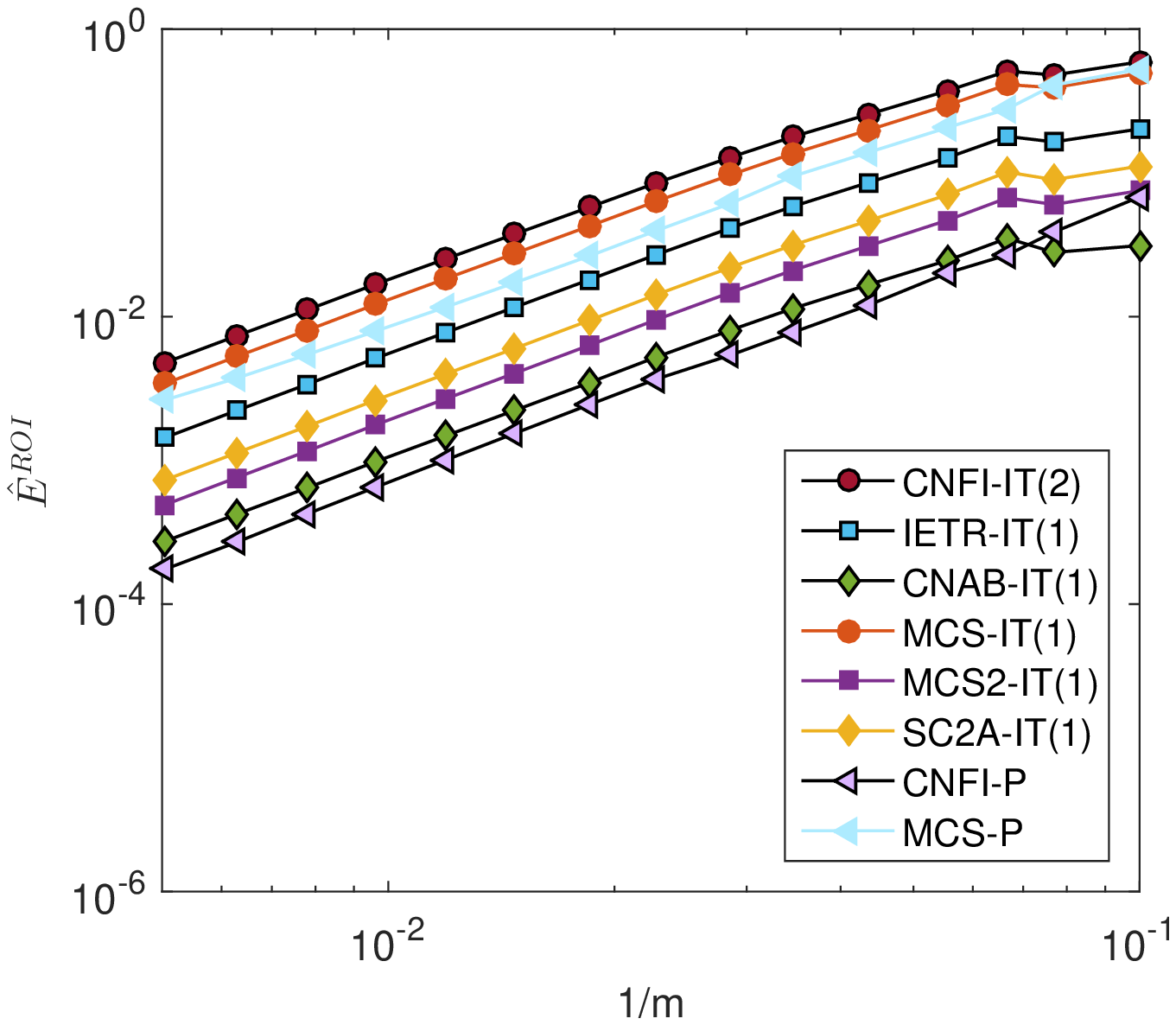}
    \includegraphics[trim={0cm 0cm 0cm 0cm},clip,
    width =0.51\textwidth]{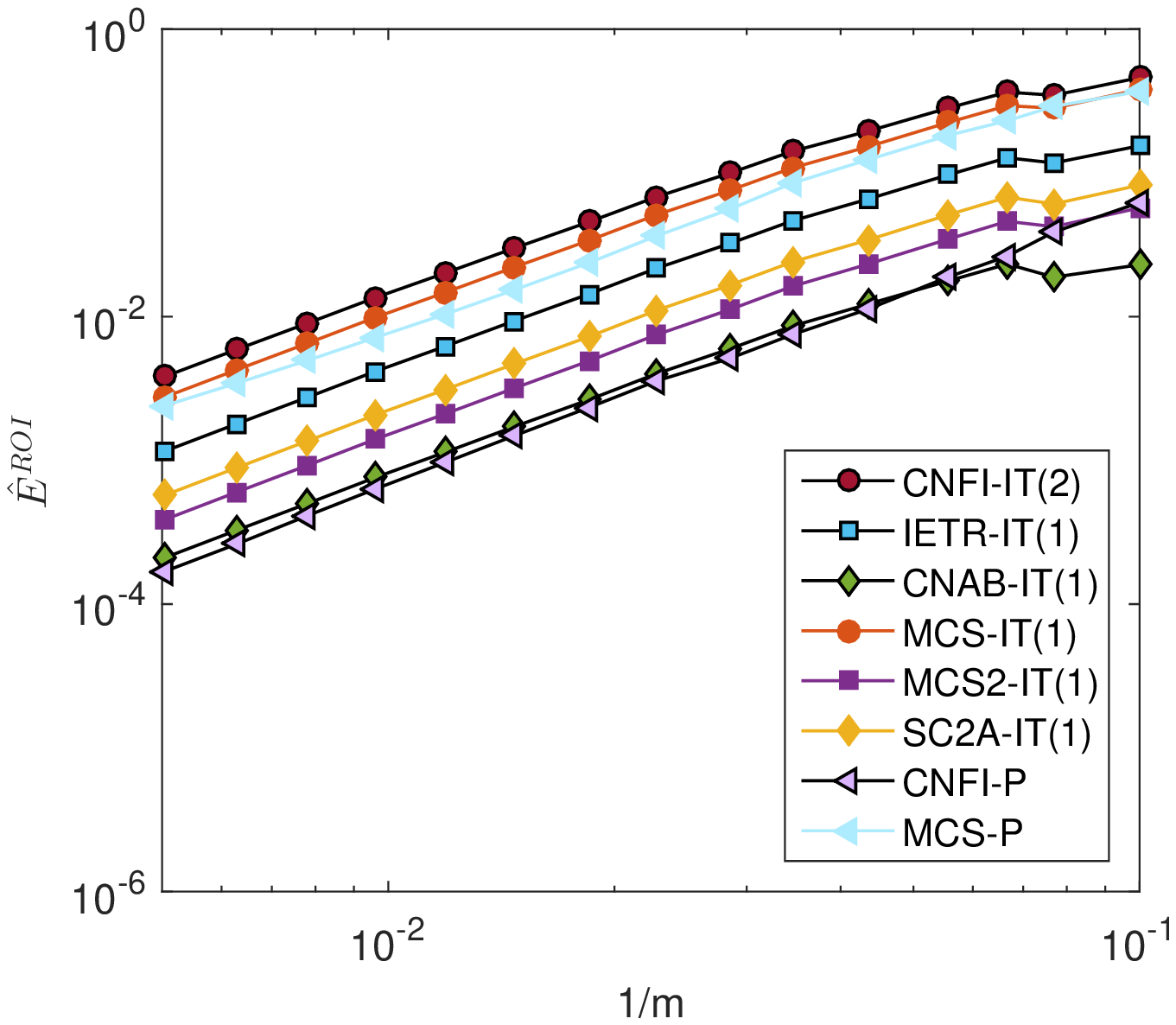}\\
    \caption{Temporal errors \eqref{temperror} of the eight operator splitting methods in the case of the American put-on-the-min option 
    under the two-asset Merton jump-diffusion model with $\kappa = 2$ for method \eqref{CNFI IT} and $\kappa = 1$ for methods 
    \eqref{IETR IT}-\eqref{SC2A IT}.
    Displayed are the errors on the large ROI (left) and the small ROI (right) and for parameter Set~1 (top), Set~2 (mid) and Set~3 (bottom) 
    from Table~\ref{tabpars}.}
    \label{fig:temperror_one_it_pom}
\end{figure}

\begin{figure}[h!]
    \centering
    \hspace{-0.5cm}\includegraphics[trim={0cm 0cm 0cm 0cm},clip,
    width =0.51\textwidth]{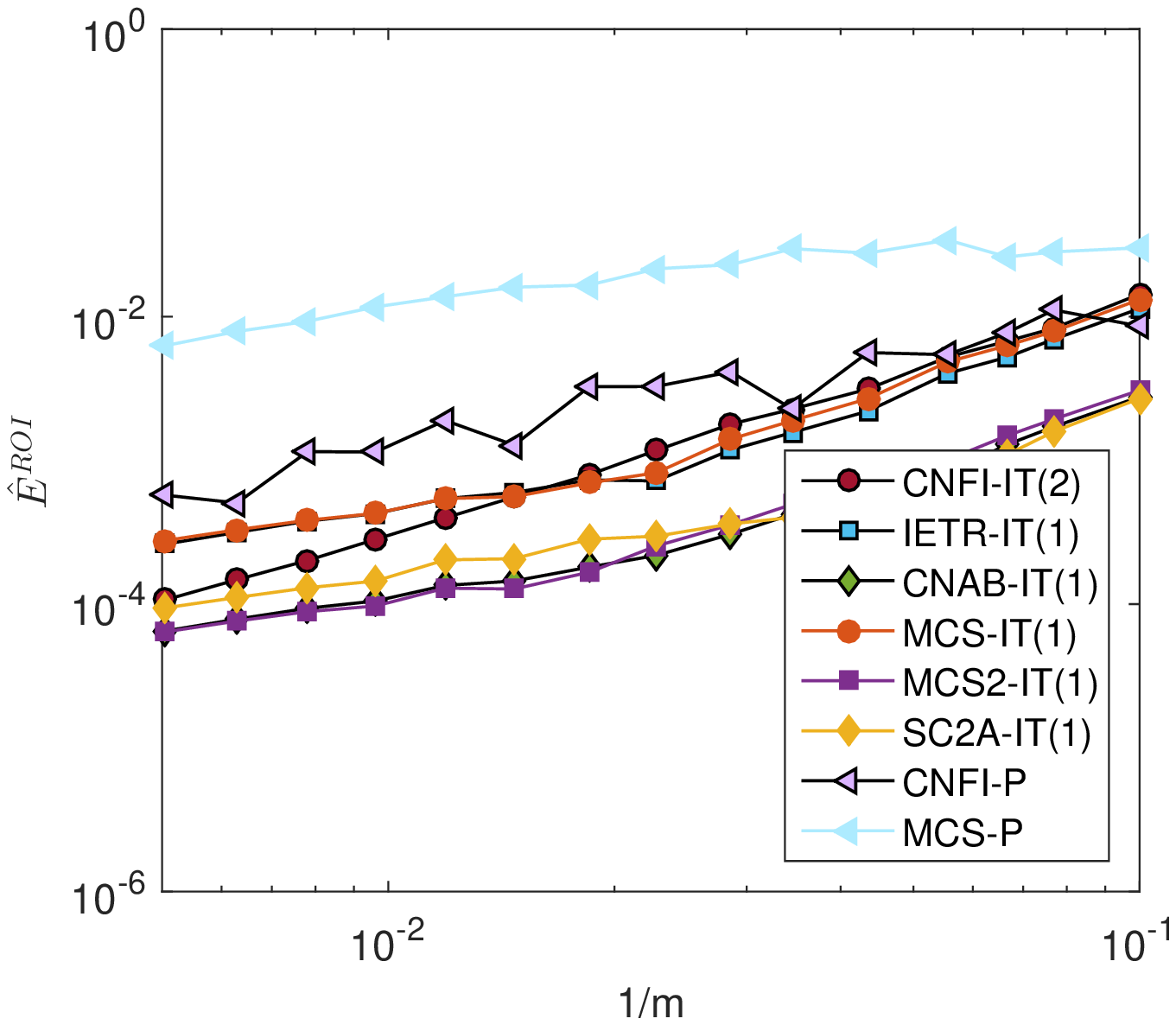}
    \includegraphics[trim={0cm 0cm 0cm 0cm},clip,
    width =0.51\textwidth]{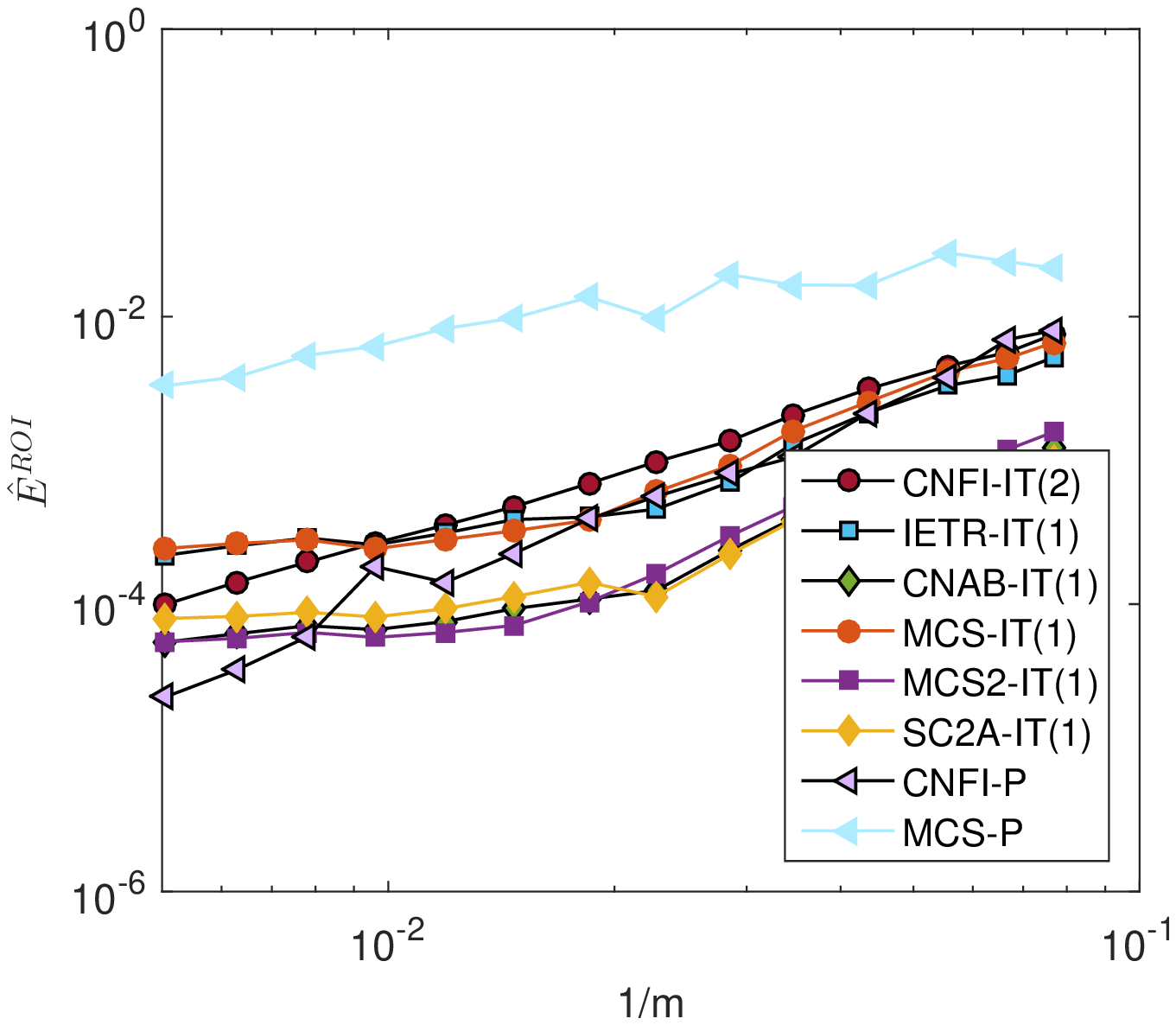}\\
    \hspace{-0.5cm}\includegraphics[trim={0cm 0cm 0cm 0cm},clip,
    width =0.51\textwidth]{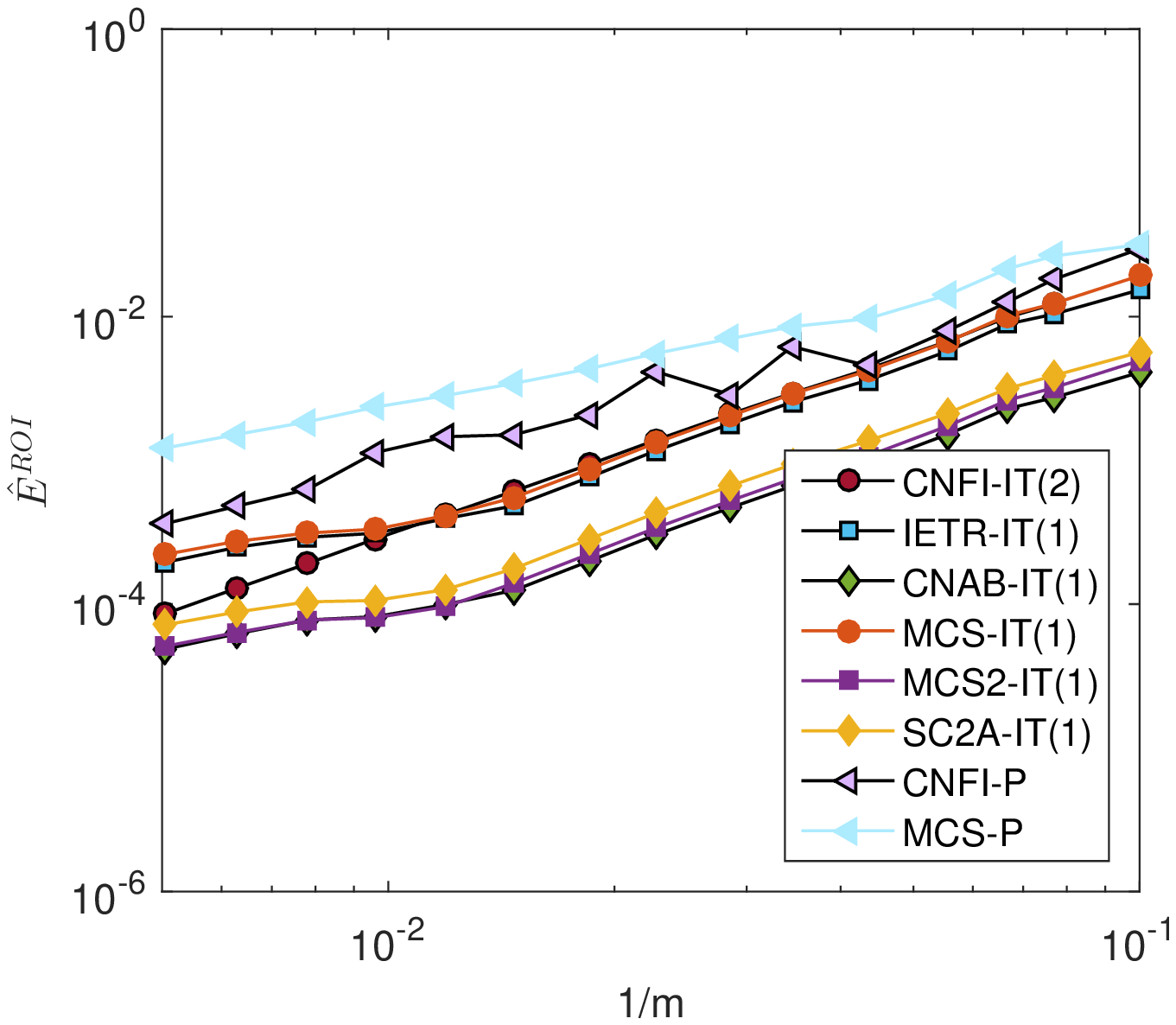}
    \includegraphics[trim={0cm 0cm 0cm 0cm},clip,
    width =0.51\textwidth]{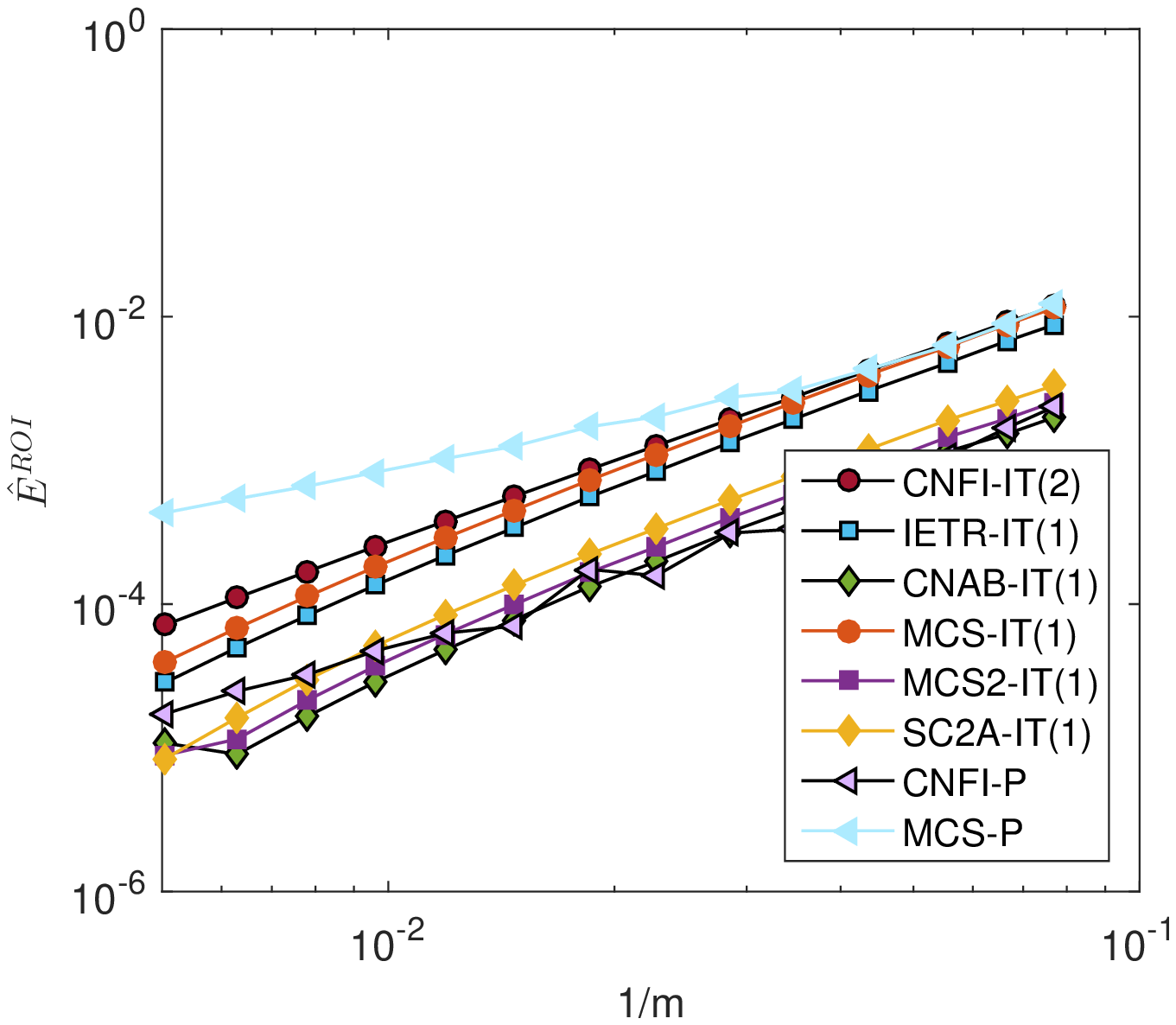}\\
    \hspace{-0.5cm}\includegraphics[trim={0cm 0cm 0cm 0cm},clip,
    width =0.51\textwidth]{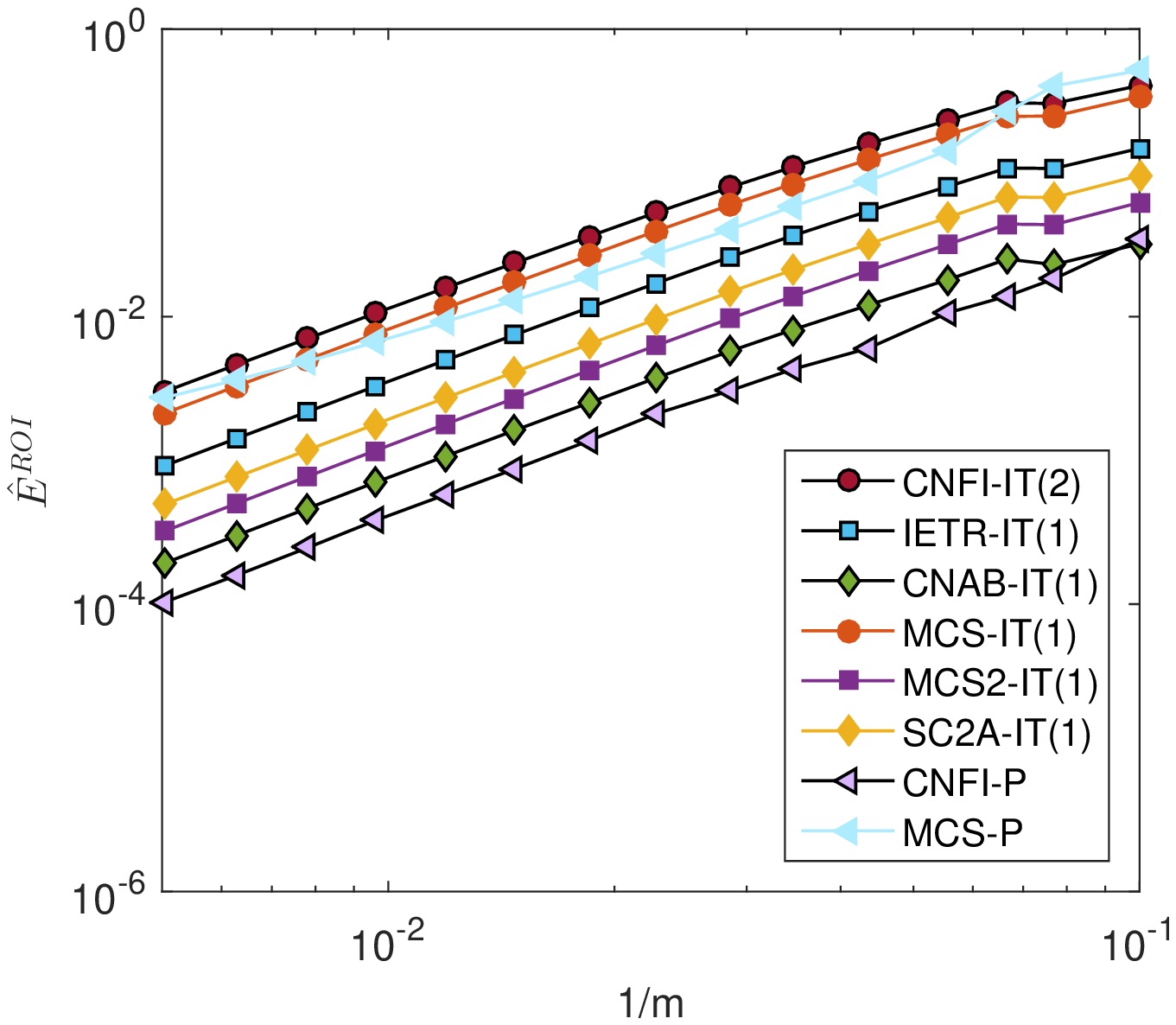}
    \includegraphics[trim={0cm 0cm 0cm 0cm},clip,
    width =0.51\textwidth]{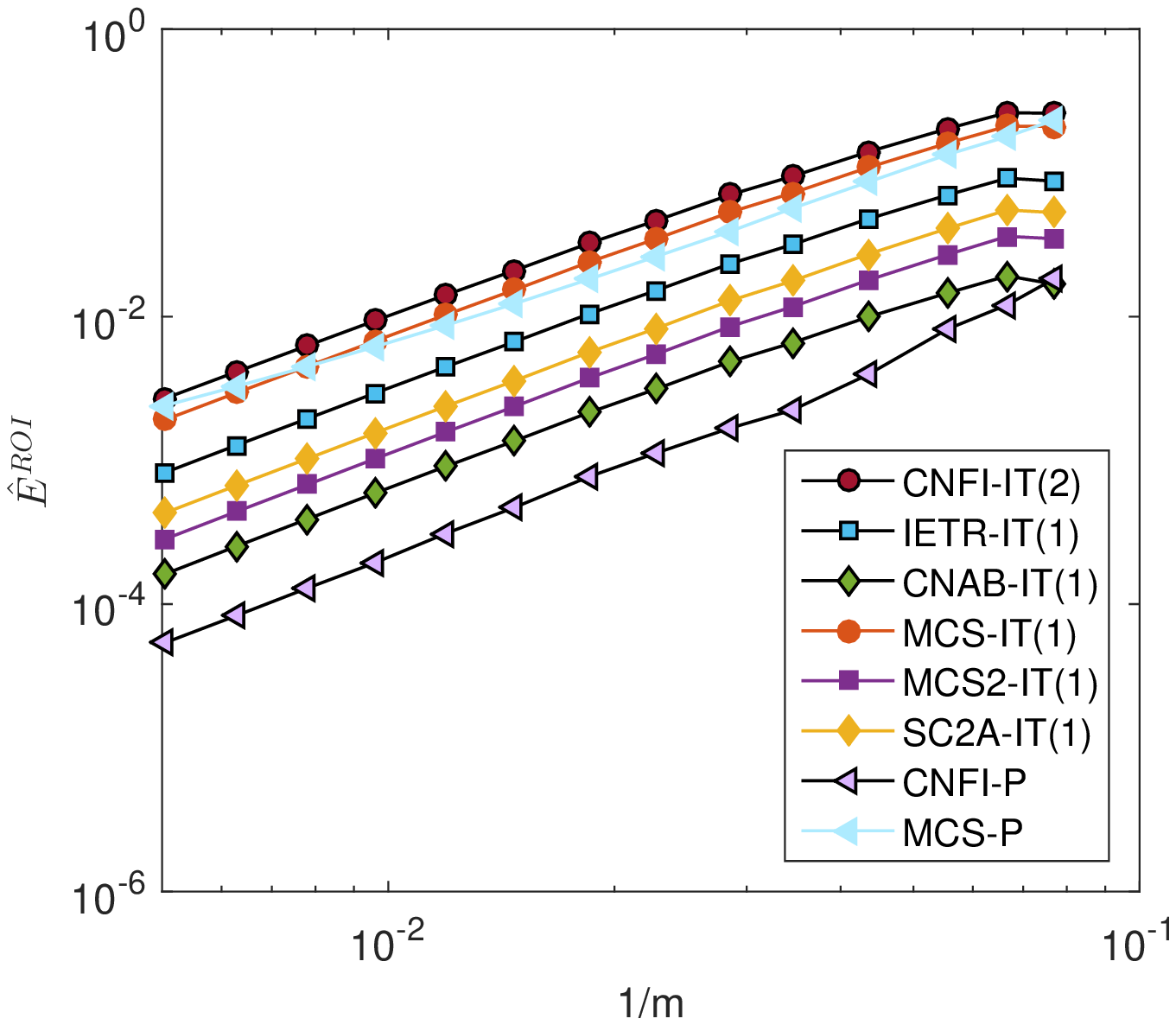}\\
    \caption{Temporal errors \eqref{temperror} of the eight operator splitting methods in the case of the American put-on-the-average option 
    under the two-asset Merton jump-diffusion model with $\kappa = 2$ for method \eqref{CNFI IT} and $\kappa = 1$ for methods 
    \eqref{IETR IT}-\eqref{SC2A IT}.
    Displayed are the errors on the large ROI (left) and the small ROI (right) and for parameter Set~1 (top), Set~2 (mid) and Set~3 (bottom) 
    from Table~\ref{tabpars}.}    
    \label{fig:temperror_one_it_poa}
\end{figure}

\begin{figure}[h!]
    \centering
    \hspace{-0.5cm}\includegraphics[trim={0cm 0cm 0cm 0cm},clip,
    width =0.51\textwidth]{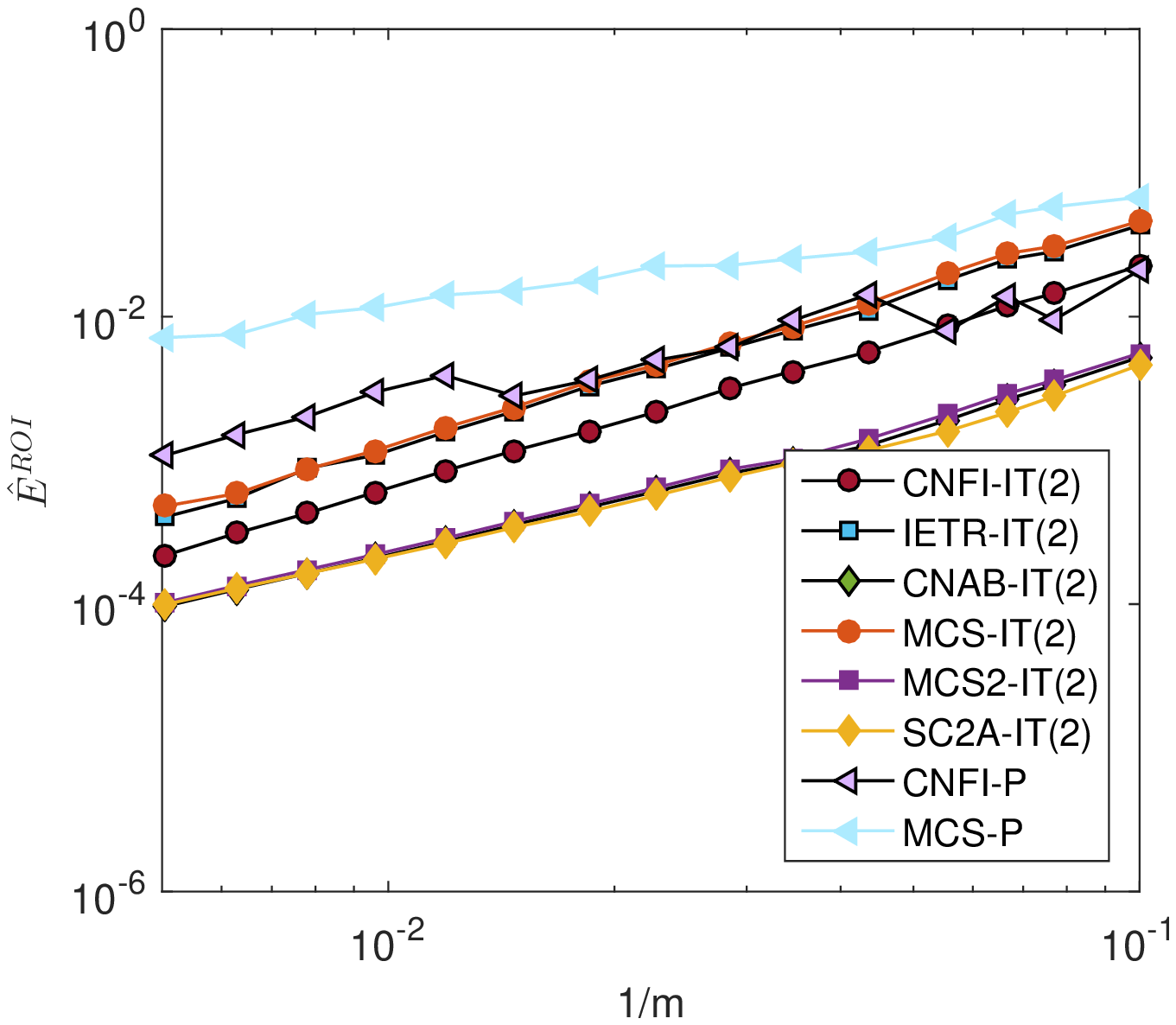}
    \includegraphics[trim={0cm 0cm 0cm 0cm},clip,
    width =0.51\textwidth]{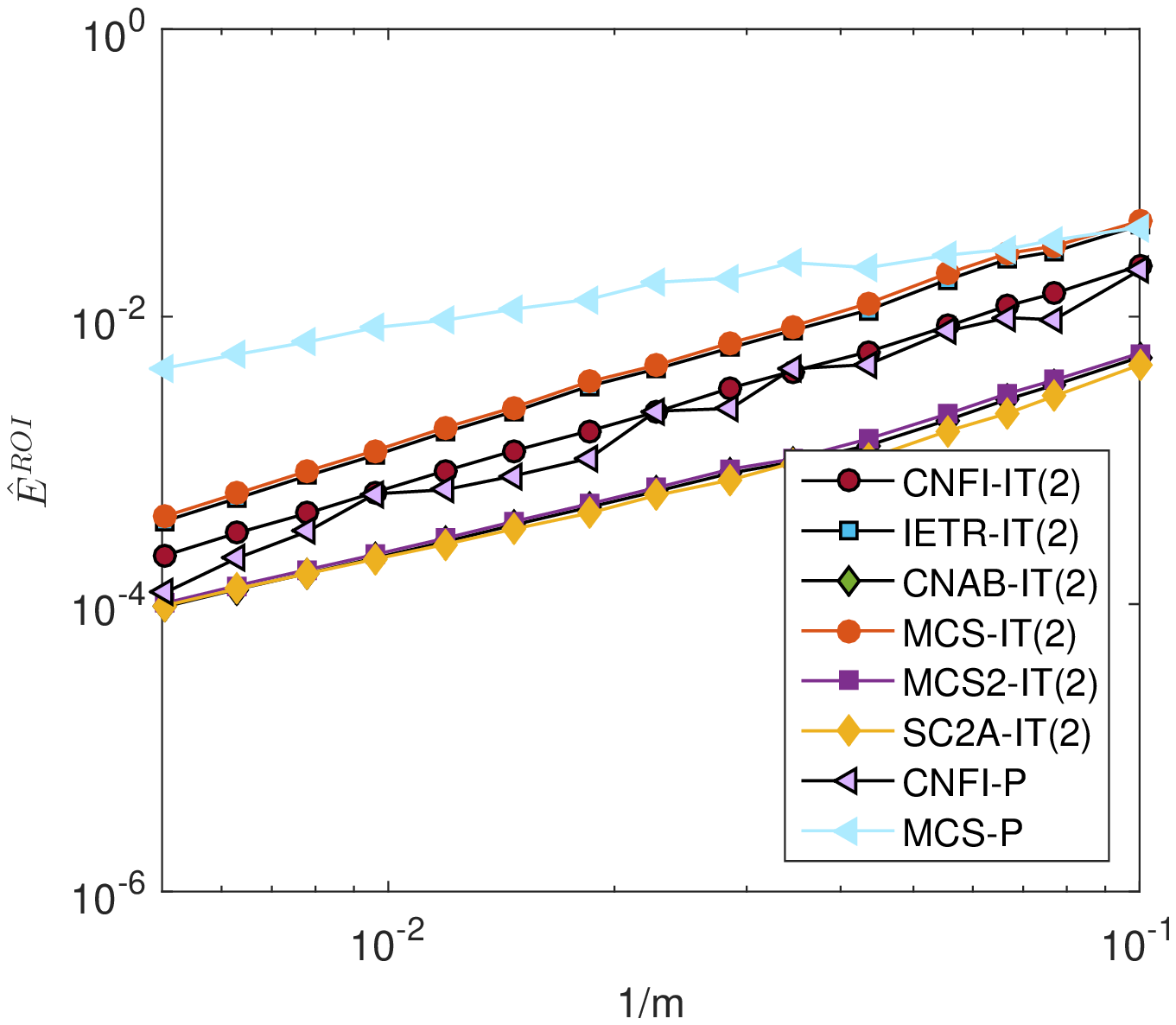}\\
    \hspace{-0.5cm}\includegraphics[trim={0cm 0cm 0cm 0cm},clip,
    width =0.51\textwidth]{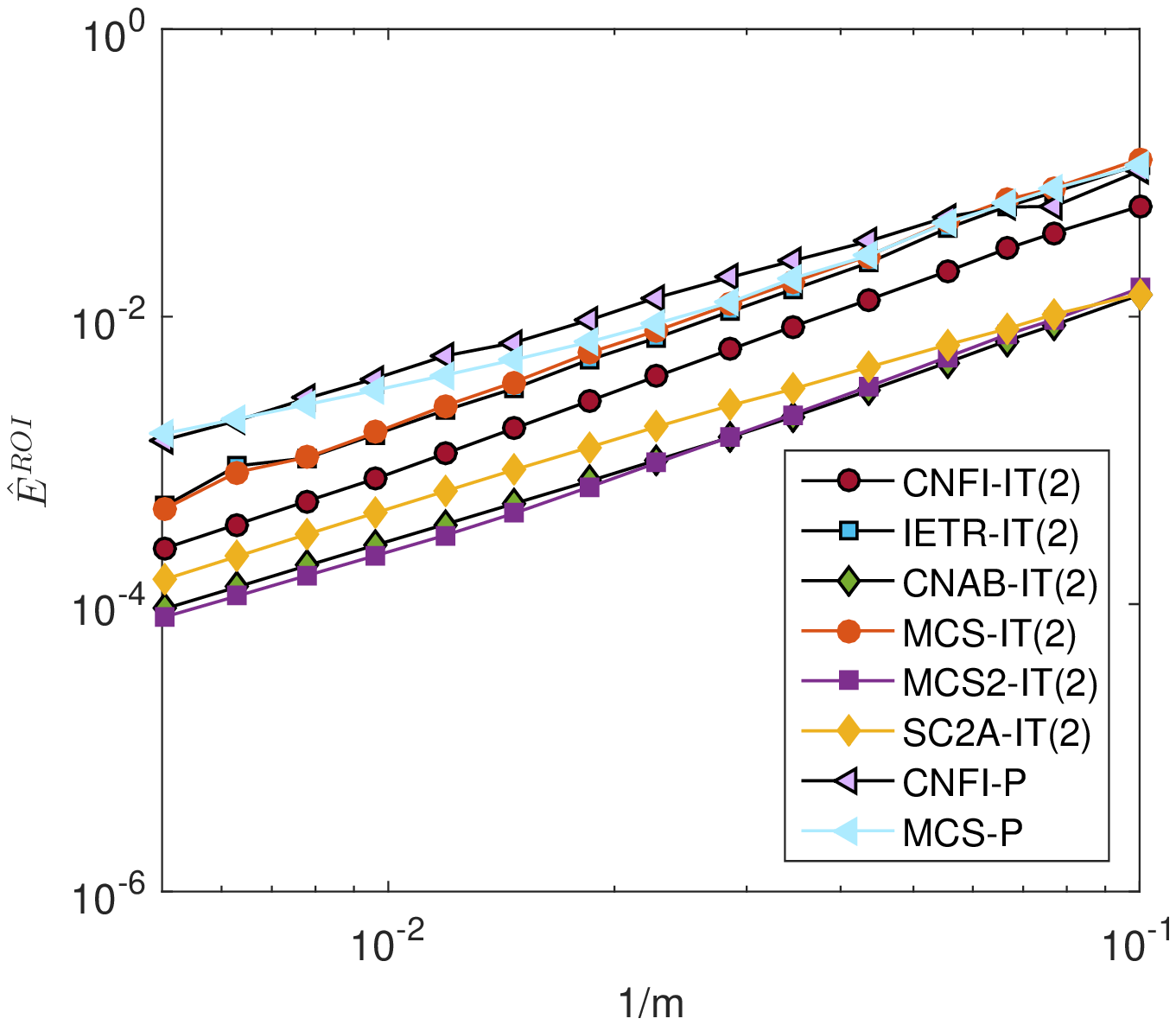}
    \includegraphics[trim={0cm 0cm 0cm 0cm},clip,
    width =0.51\textwidth]{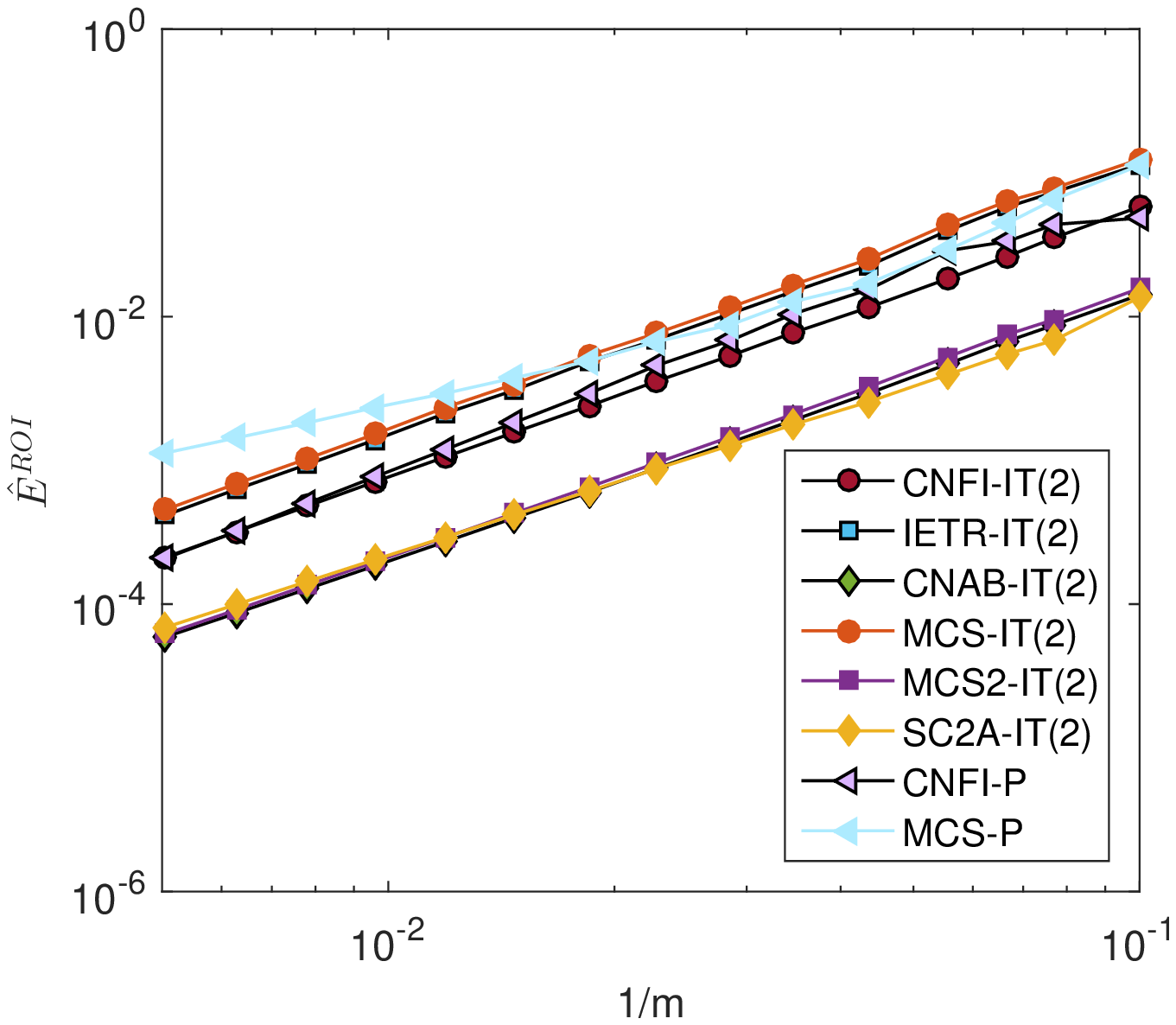}\\
    \hspace{-0.5cm}\includegraphics[trim={0cm 0cm 0cm 0cm},clip,
    width =0.51\textwidth]{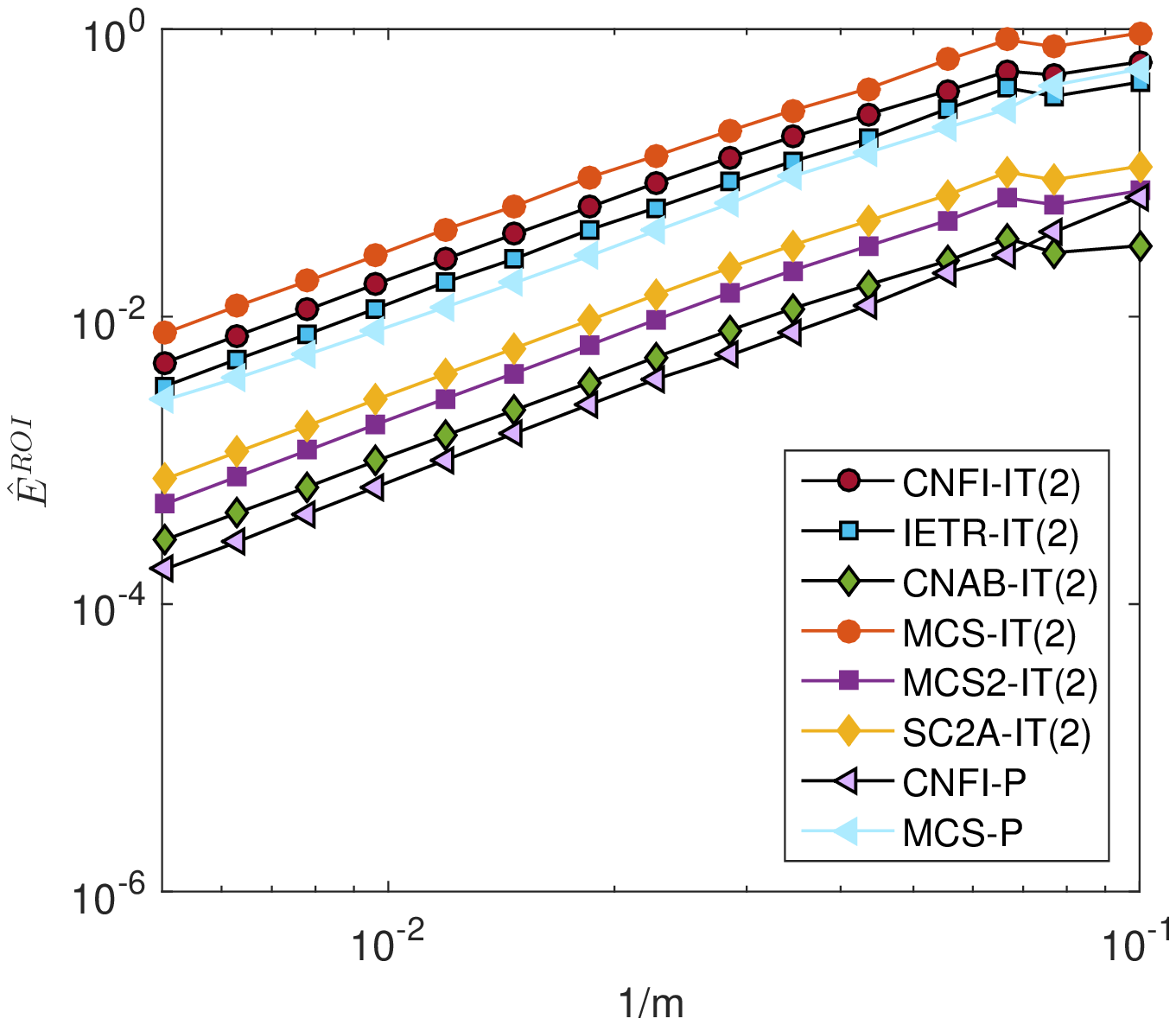}
    \includegraphics[trim={0cm 0cm 0cm 0cm},clip,
    width =0.51\textwidth]{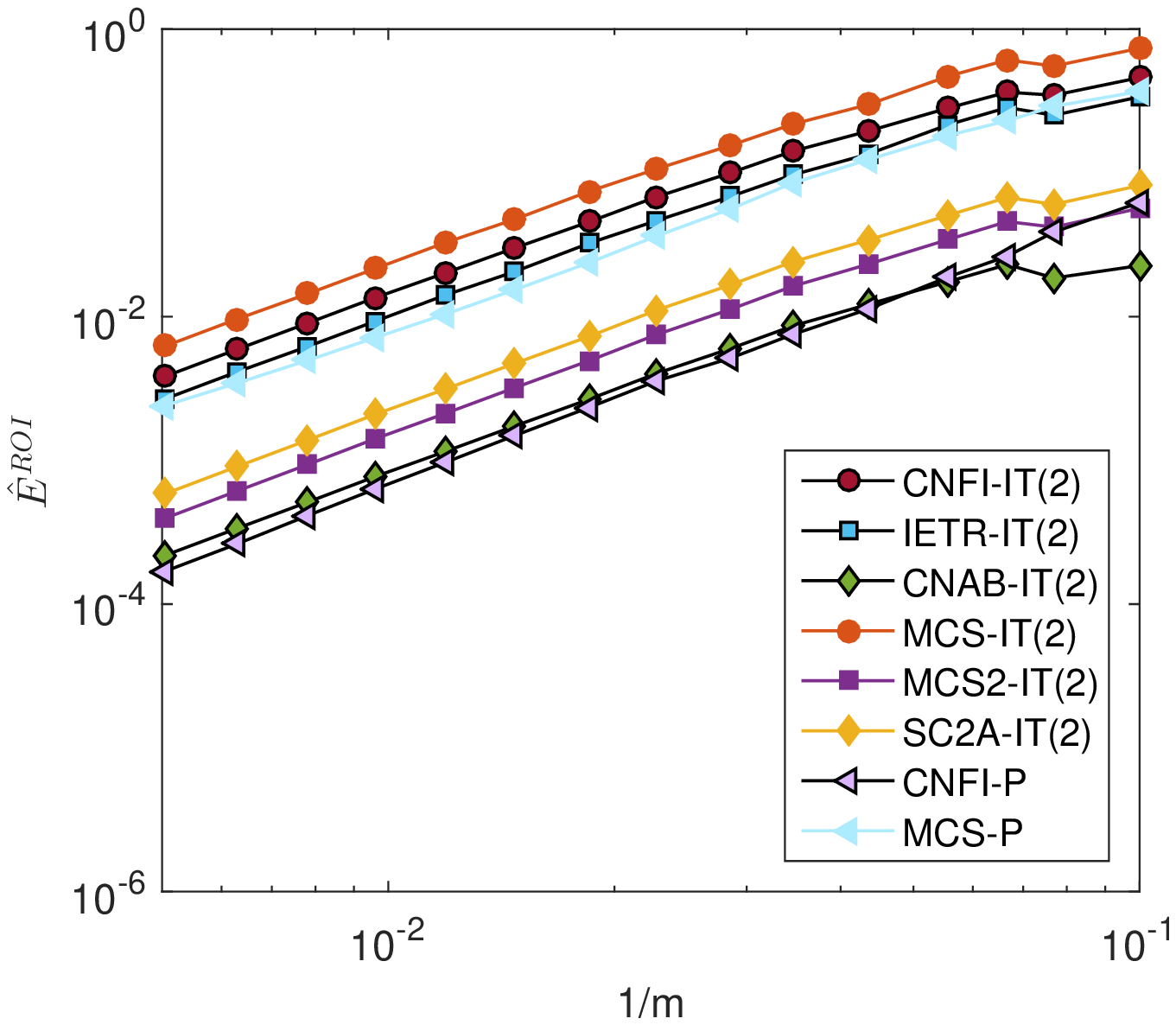}\\
    \caption{Temporal errors \eqref{temperror} of the eight operator splitting methods in the case of the American put-on-the-min option 
    under the two-asset Merton jump-diffusion model with $\kappa = 2$ for all methods \eqref{CNFI IT}-\eqref{SC2A IT}.
    Displayed are the errors on the large ROI (left) and the small ROI (right) and for parameter Set~1 (top), Set~2 (mid) and Set~3 (bottom) 
    from Table~\ref{tabpars}.}
    \label{fig:temperror_pom}
\end{figure}

\begin{figure}[h!]
    \centering
    \hspace{-0.5cm}\includegraphics[trim={0cm 0cm 0cm 0cm},clip,
    width =0.51\textwidth]{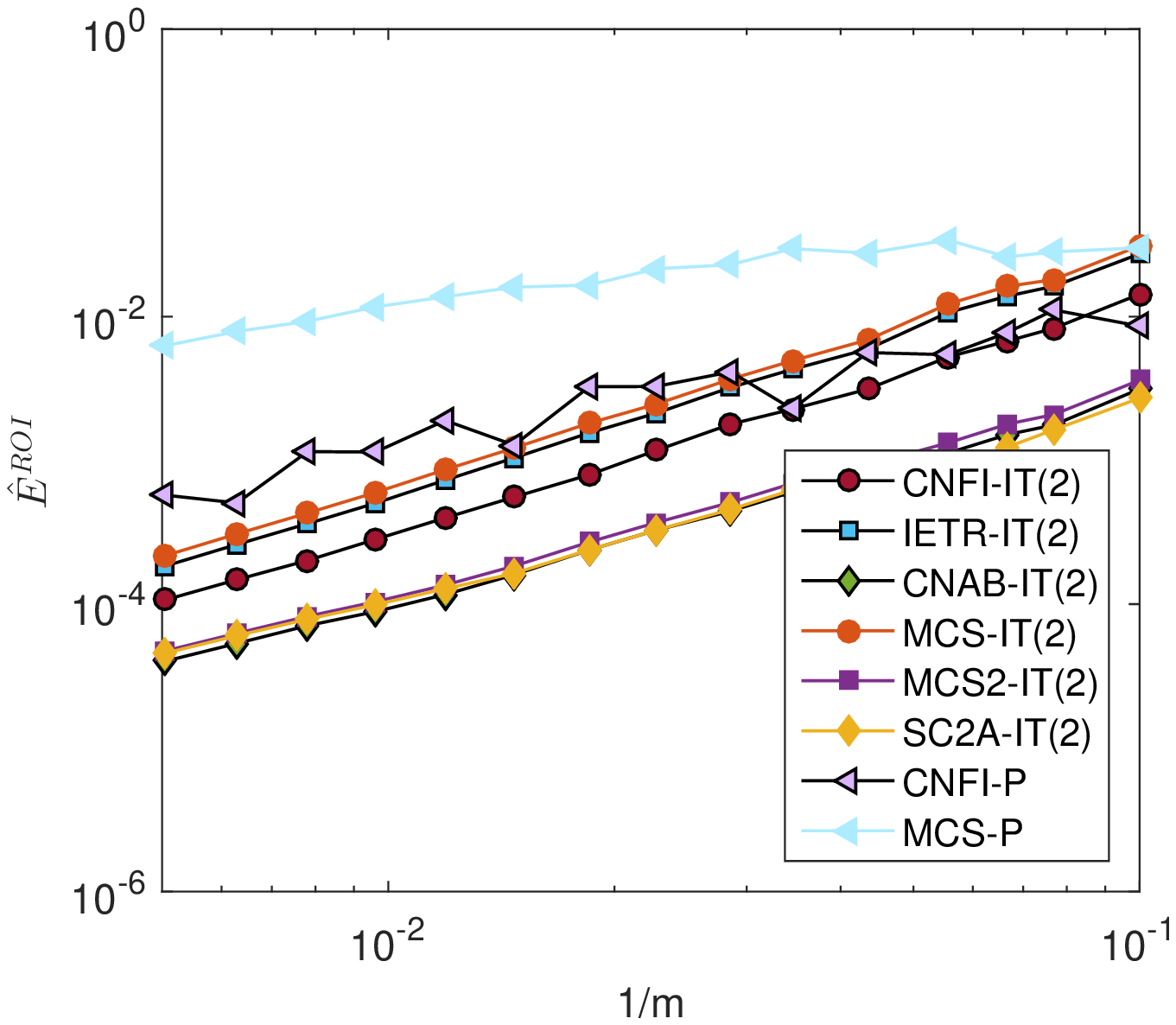}
    \includegraphics[trim={0cm 0cm 0cm 0cm},clip,
    width =0.51\textwidth]{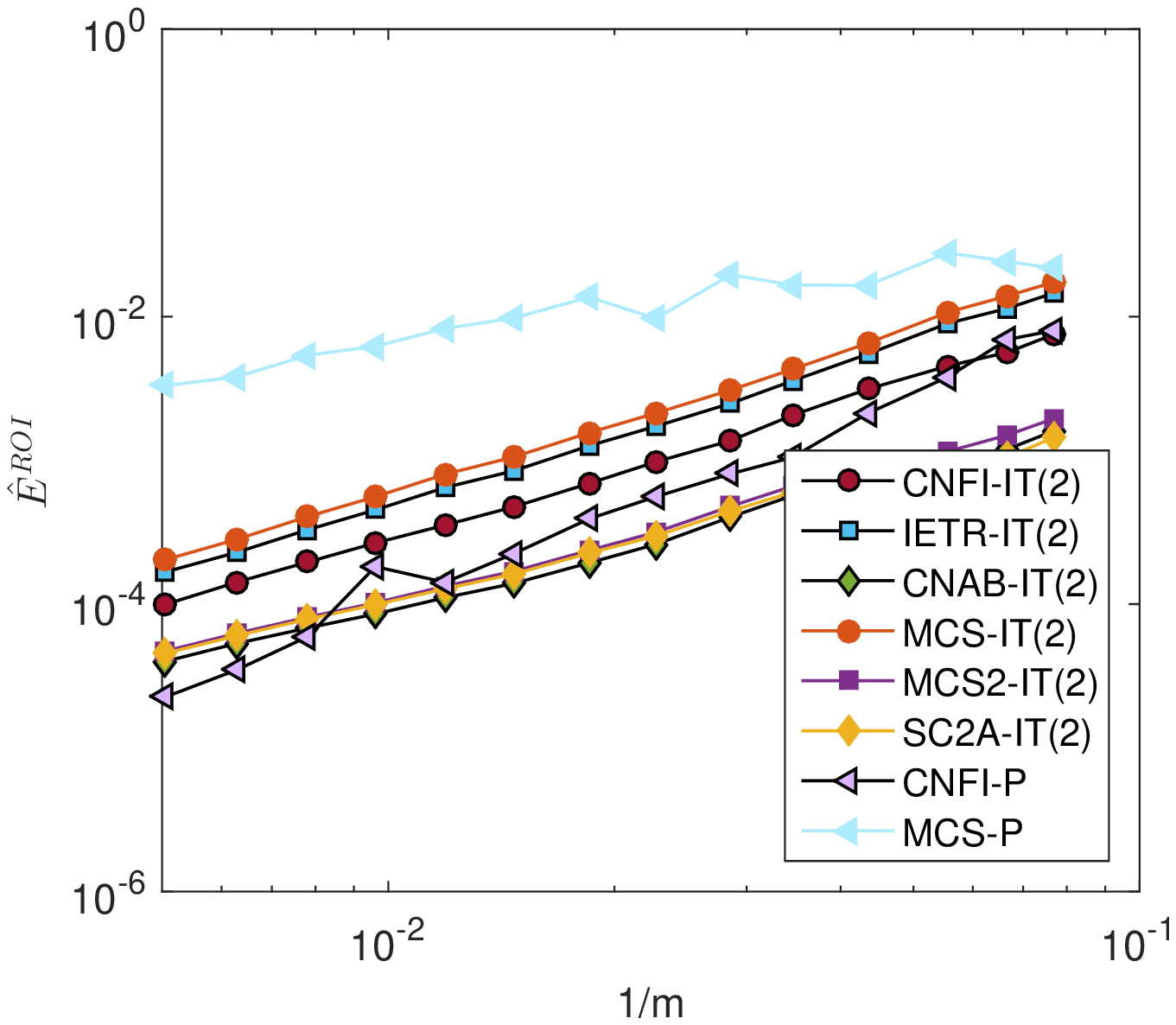}\\
    \hspace{-0.5cm}\includegraphics[trim={0cm 0cm 0cm 0cm},clip,
    width =0.51\textwidth]{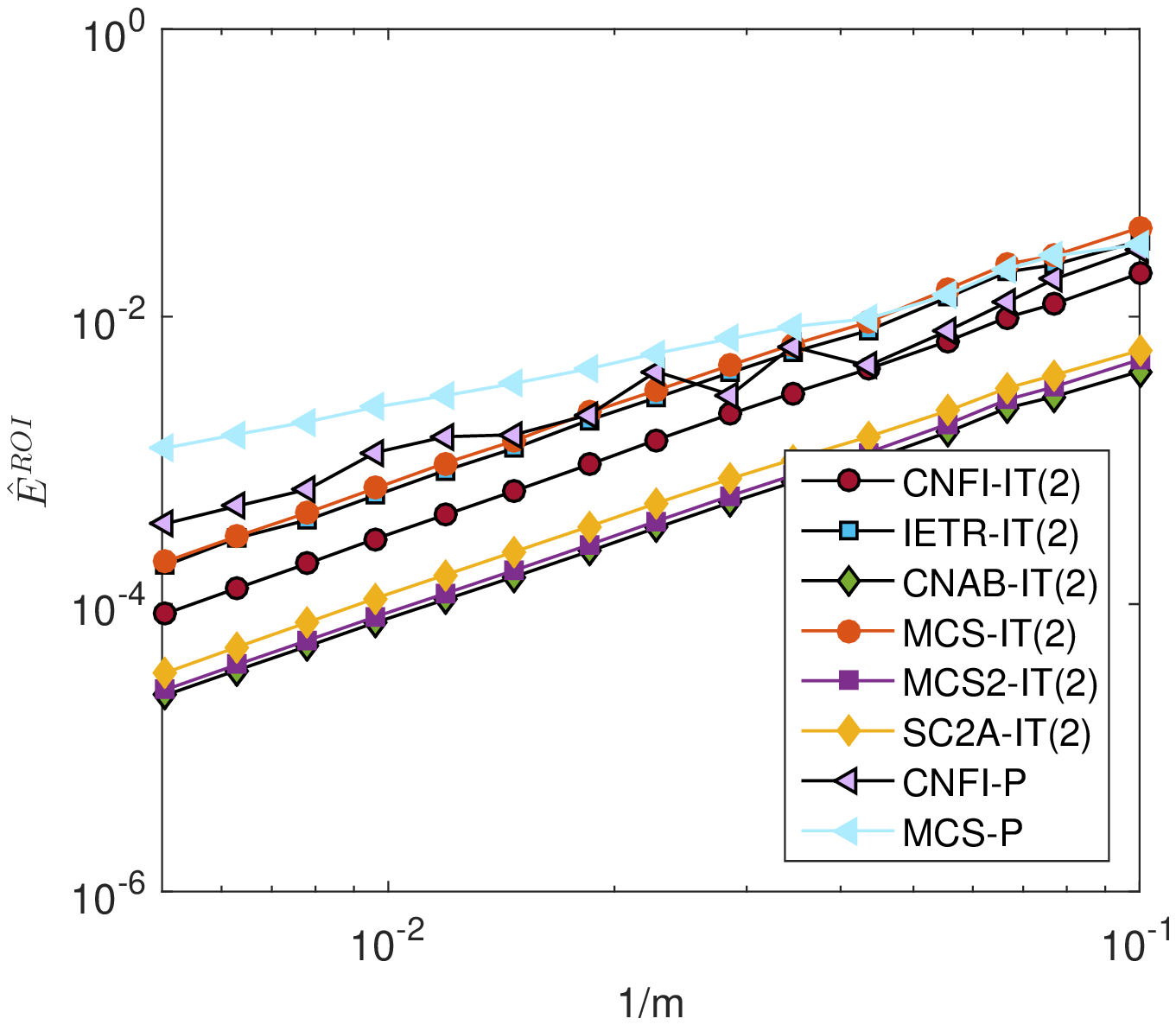}
    \includegraphics[trim={0cm 0cm 0cm 0cm},clip,
    width =0.51\textwidth]{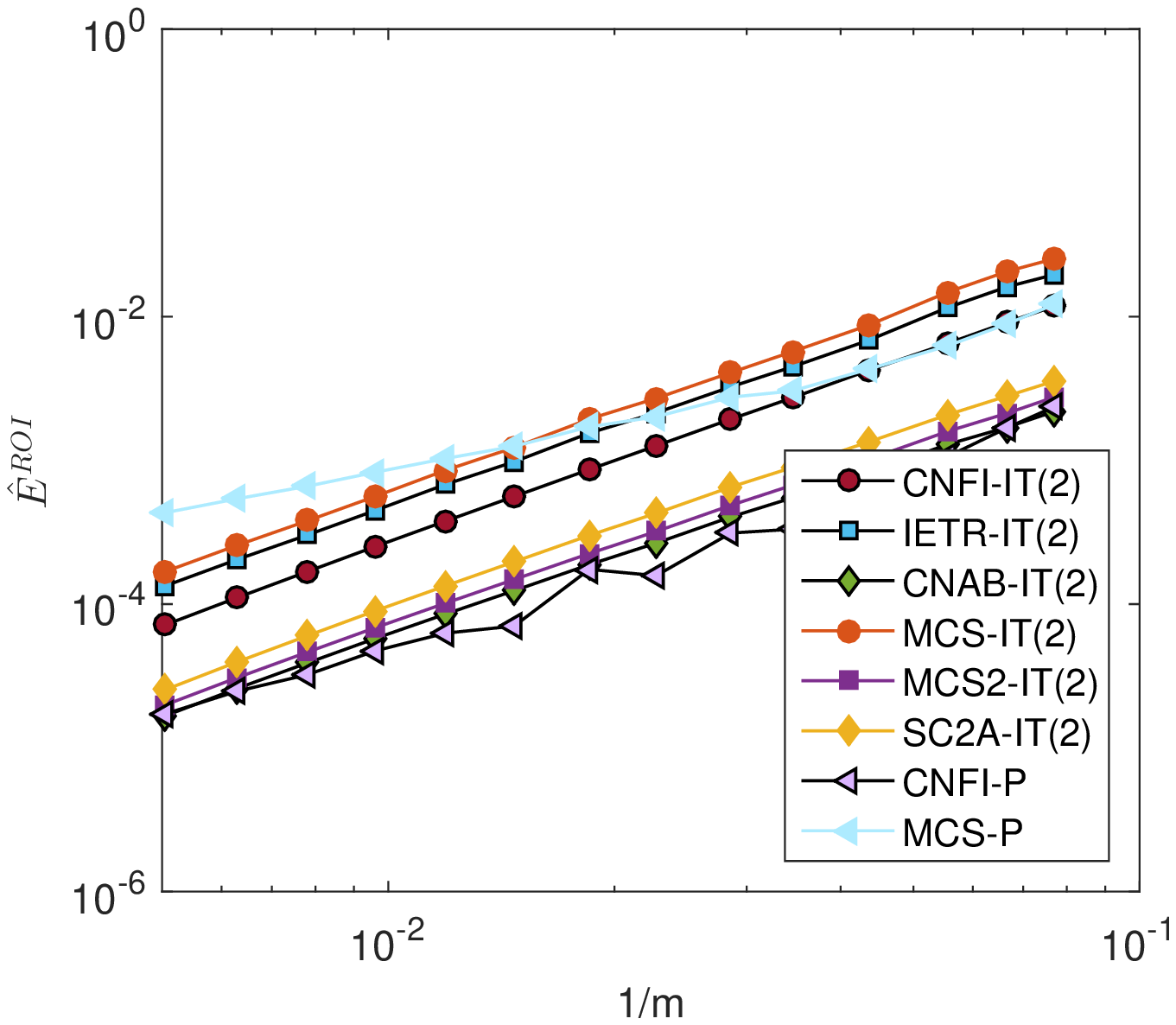}\\
    \hspace{-0.5cm}\includegraphics[trim={0cm 0cm 0cm 0cm},clip,
    width =0.51\textwidth]{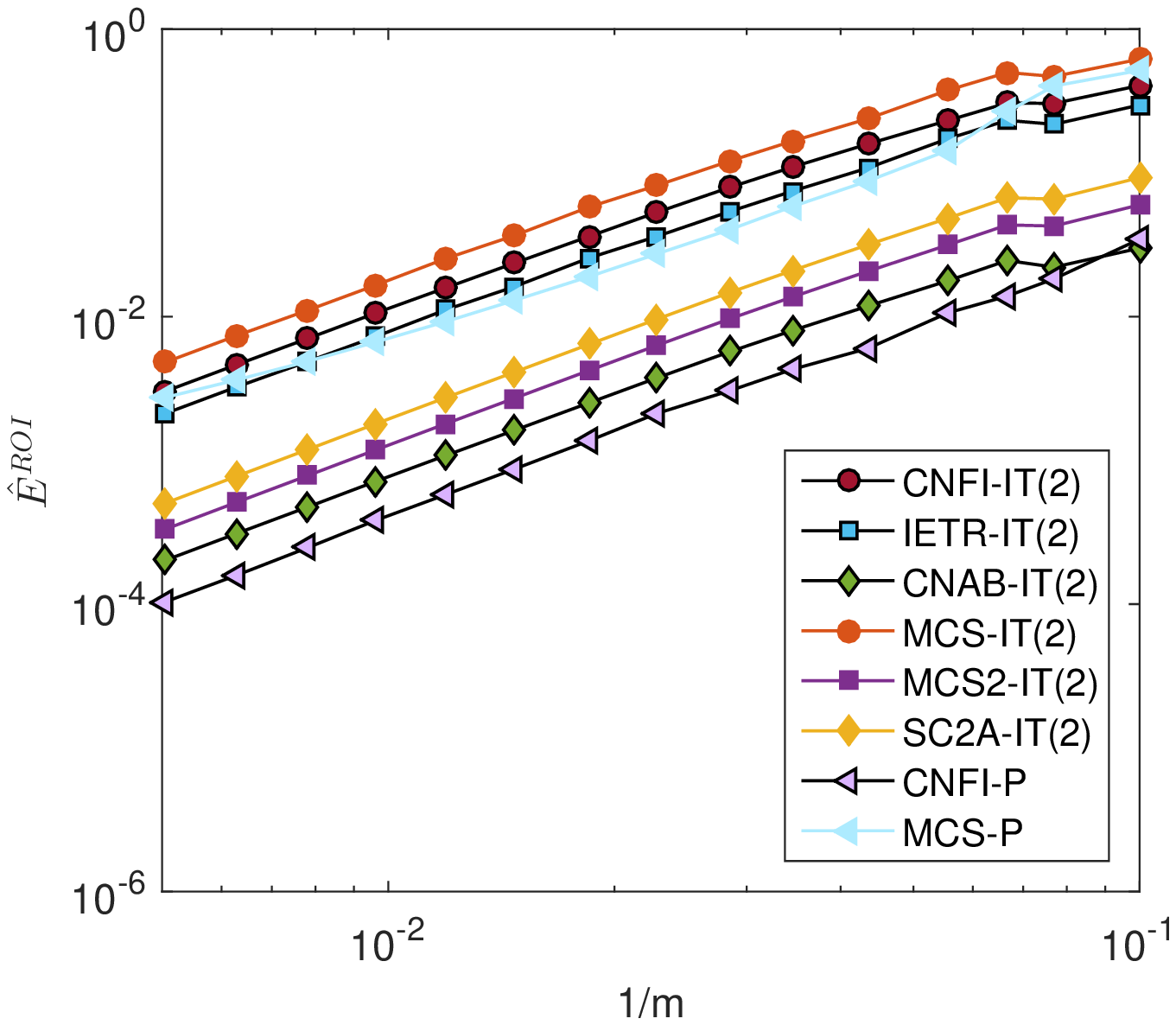}
    \includegraphics[trim={0cm 0cm 0cm 0cm},clip,
    width =0.51\textwidth]{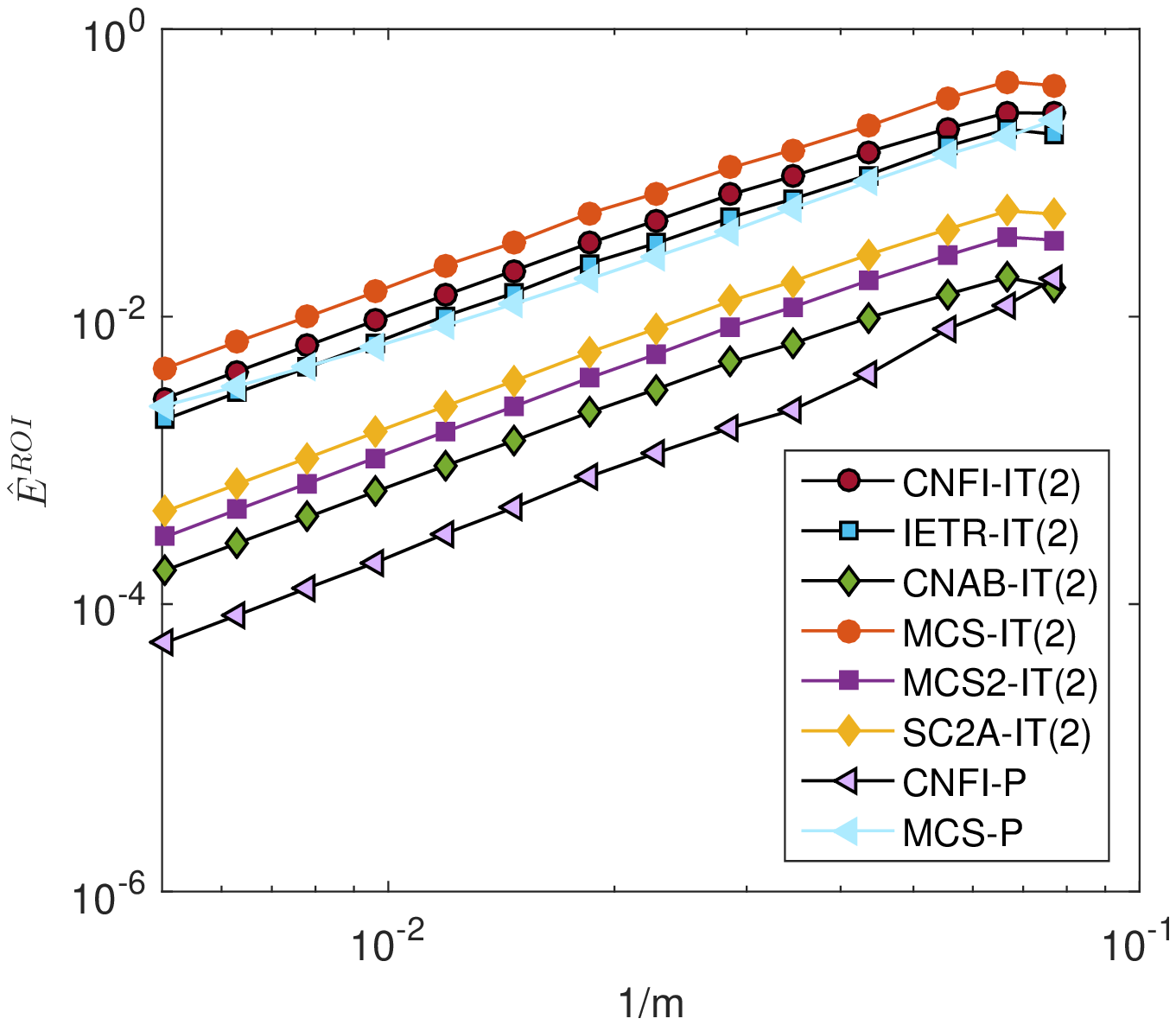}\\
    \caption{Temporal errors \eqref{temperror} of the eight operator splitting methods in the case of the American put-on-the-average option 
    under the two-asset Merton jump-diffusion model with $\kappa = 2$ for all methods \eqref{CNFI IT}-\eqref{SC2A IT}.
    Displayed are the errors on the large ROI (left) and the small ROI (right) and for parameter Set~1 (top), Set~2 (mid) and Set~3 (bottom) 
    from Table~\ref{tabpars}.}
    \label{fig:temperror_poa}
\end{figure}

In this section, ample numerical experiments are performed to gain insight into the convergence behaviour of the eight operator splitting 
methods formulated in Section~\ref{tempdisc} in the numerical solution of the semidiscrete two-dimensional Merton PIDCP \eqref{PIDCdisc}. 
To this purpose, the {\it temporal discretization error} is considered at $t=T$ on a region of interest (ROI) in the spatial domain,
\begin{align}\label{temperror}
\widehat{E}^{ROI}(m,N) = \max \left\{|\wV^{N'}_{i,j}-V_{i,j}(T)|: \Delta t = T/N'~\textrm{and}~s_L < s_{1,i},s_{2,j}< s_U \right\}.
\end{align}
Here $s_L$ and $s_U$ denote given, fixed numbers with $0< s_L < s_U$ that define the ROI.
We take $m_1=m_2=m$ for the spatial discretization.
The number of time steps $N'$ is chosen in function of $N$ and the splitting method so that the total number of matrix-vector products 
with the matrix $A^{(J)}$ over the whole time interval $[0,T]$ is (essentially) the same for all methods.
This is done to have a fair comparison between the eight methods, as determining these matrix-vector products forms the dominating 
computational part of each time step.
Accordingly, we take $N' = \lceil 2N/\kappa \rceil$ for method \eqref{CNFI IT}, $N'= \lceil 2N/(\kappa+1)\rceil$ for methods \eqref{IETR IT} 
and \eqref{MCS IT}, $N'= 2N$ for \eqref{CNAB IT}, \eqref{MCS2 IT} and \eqref{SC2A IT} and $N' = N$ for \eqref{MCS P}.
For the CNFI-P method \eqref{theta-P} the number of matrix-vector products with $A^{(J)}$ per time step is not known a priori, since a 
dynamic convergence criterion \eqref{convcrit} is employed. 
However, numerical experiments reveal that convergence of the penalty iteration is attained after on average four iterations. 
In view of this, we consider method \eqref{theta-P} with $N' = \lceil N/2 \rceil$. 
A reference solution $V(T)$ to the semidiscrete PIDCP \eqref{PIDCdisc} has been computed by applying the CNFI-P method\footnote{Using 
e.g.~the CNAB-IT(2) method leads to visually identical Figures \ref{EERs}--\ref{fig:temperror_poa}.} and $N' = 10N$ time steps.

We deal with two types of options: an American put-on-the-min option and an American put-on-the-average option. 
Their payoff functions, with given strike $K$, are 
\begin{align*}
    \phi_{\text{put-on-min}}(s_1,s_2) = \max(0\,,\,K-\min(s_1,s_2))
\end{align*}
and
\begin{align*}
    \phi_{\text{put-on-average}}(s_1,s_2) = \max\left(0\,,\,K-\frac{s_1+s_2}{2}\right).
\end{align*}
Three financial parameter sets are considered, which are specified in Table~\ref{tabpars}. 
They are identical to those chosen in Boen \& in 't Hout~\cite{BH19} for the case of European options.
The first parameter set is given in Clift \& Forsyth~\cite{CF08}.
Here $\lambda T$ is small, indicating a low expected number of jumps. 
The second set has the same diffusion parameters as in Zvan, Forsyth \& Vetzal~\cite{ZFV01} and jump parameters are taken 
where $\lambda T$ is about the same size as for the first set.
The third set was introduced in~\cite{BH19} and here $\lambda T$ is quite large.
Note further that for all three sets the correlation coefficients $\rho$ and $\rhoh$ are nonzero.

\begin{table}[h!]
\centering
\begin{tabular}{ c | c   c  c c c c c c c c c c }
 & $\sigma_1$ & $\sigma_2$ & $\rho$ & $\lambda$ & $\gamma_1$ & $\gamma_2$ & $\rhoh$ & $\delta_1$ & $\delta_2$ & $r$ & $K$ & $T$ \\
 \hline
Set~1 & 0.12 & 0.15 & 0.30 & 0.60 & -0.10 & 0.10 & -0.20 & 0.17 & 0.13 & 0.05 & 100 & 1\\
Set~2 & 0.30 & 0.30 & 0.50 & 2 & -0.50 & 0.30 & -0.60 & 0.40 & 0.10 & 0.05 & 40 & 0.5\\
Set~3 & 0.20 & 0.30 & 0.70 & 8 & -0.05 & -0.20 & 0.50 & 0.45 & 0.06 & 0.05 & 40 & 1
 \end{tabular}
\caption{Parameter sets for the two-asset Merton jump-diffusion model and American option.}\label{tabpars}
\end{table}

Figure~\ref{EERs} displays in grey the approximated early exercise regions (EERs) for the American put-on-the-min and put-on-the-average options 
under the three parameter sets from Table~\ref{tabpars}.
In the following we study, for all eight methods formulated in Section~\ref{tempdisc}, the temporal discretization error \eqref{temperror} on 
both a large and a small ROI.
These regions are indicated in Figure~\ref{EERs} with blue and red, respectively.
The large ROI is given by $[(1/2)K,(3/2)K]$ and has a nonempty intersection with the EERs for parameter Sets~1 and 2.
The small ROI is given by $[(7/8)K,(9/8)K]$ in the case of the put-on-the-min option and $[(9/10)K,(11/10)K]$ in the case of the put-on-the-average 
option.
From Figure~\ref{EERs} we observe that the small ROI does not intersect the EERs, that is, it lies fully within the continuation regions.

We start by considering the methods \eqref{IETR IT}-\eqref{SC2A IT} with IT$(\kappa)$ splitting and $\kappa = 1$.
Only the CNFI-IT($\kappa$) method \eqref{CNFI IT} is used with $\kappa = 2$, so that the order of consistency of the underlying IMEX scheme is also 
equal to two, compare Subsection~\ref{IT schemes}.
As discussed above, for a fair comparison, the methods \eqref{CNFI IT}, \eqref{IETR IT}, \eqref{MCS IT} are applied in this case with $N'=N$ time 
steps and the methods \eqref{CNAB IT}, \eqref{MCS2 IT}, \eqref{SC2A IT} with $N' = 2N$ time steps. 
We take $N=m$ and consider a range of values $m$ with $10\leq m \leq 200$.
The temporal errors of all eight methods are shown (versus $1/m$) in Figure~\ref{fig:temperror_one_it_pom} for the put-on-the-min option and 
in Figure~\ref{fig:temperror_one_it_poa} for the put-on-the-average option, under the three parameter sets given in Table~\ref{tabpars}. 
Here the left column corresponds to the large ROI and the right column the small ROI.
As a positive observation, all eight temporal discretization methods show a regular, monotonic convergence behaviour for parameter Set~3.
However, with the favourable exception of the CNFI-IT(2) method, each of the other methods can yield a less regular convergence behaviour 
under Sets~1 and 2. In particular, the temporal errors can level off as $m$ increases. 
This is found notably in the cases where the ROI overlaps with, or lies close to, the EER.

Since the CNFI-IT$(2)$ method appears to be more robust, we next consider the methods \eqref{IETR IT}-\eqref{SC2A IT} with IT$(\kappa)$ splitting 
for $\kappa = 2$ as well.
For a fair comparison, \eqref{IETR IT} and \eqref{MCS IT} are now applied with $N' = \lceil 2N/3\rceil$ time steps (and the others as above).
In Figures~\ref{fig:temperror_pom} and \ref{fig:temperror_poa} the temporal errors are displayed for, respectively, the put-on-the-min option 
and the put-on-the-average option, under the three parameter sets from Table~\ref{tabpars} and for both the large and small ROI.
The positive result is clearly observed that if $\kappa = 2$, then each of the six methods \eqref{CNFI IT}-\eqref{SC2A IT} always shows a regular, 
monotonic convergence behaviour. 
Tables~\ref{OoC1}, \ref{OoC2}, \ref{OoC3} display the numerical orders of convergence of all methods for Sets~1, 2, 3, respectively, on the
small ROI.
They have been computed based on the ten largest values of $m$ under consideration.
For the methods \eqref{CNFI IT}-\eqref{SC2A IT} with IT$(2)$ splitting, an order of convergence is obtained which lies between 1.2 and 1.6 for 
Set~1, between 1.7 and 1.9 for Set~2, and is about equal to 1.9 for Set~3.
We conjecture that the relatively lower orders of convergence for Set~1 are related to the small ROI lying close by the EER, with the premise 
that at the boundary of the EER the option value function is less smooth.
The CNAB-IT$(2)$ and MCS2-IT$(2)$ methods, which both treat the integral part in a two-step Adams--Bashforth fashion, always gave 
rise to the smallest error constant in our experiments among all considered methods with IT$(2)$ splitting.

We mention that additional experiments have been performed with an extra iteration in the methods \eqref{CNFI IT}-\eqref{SC2A IT}, that is 
$\kappa = 3$, but this did not lead to a further significant improvement of the temporal convergence behaviour.

For the CNFI-P method \eqref{theta-P}, as already alluded to above, the convergence behaviour can be less regular for Sets~1 and 2.
For each given $N=m$, the temporal error of this method is generally found to lie in, or near to, the range of those obtained with 
all IT$(2)$ splitting type methods together.

The MCS-P method \eqref{MCS P} often yields relatively large temporal errors, notably for Set~1.
We also observed this for the values of $\theta$, $\textit{tol}$, $\textit{Large}$ from \cite{HC18}, for the nonuniform temporal 
grid \eqref{nonunitemp} and when applied to the two-dimensional Black--Scholes PDCP (without integral part).

For future reference, approximations to the values of the American put-on-the-min option, respectively  put-on-the-average option,
under the three parameter sets are summarized in Table~\ref{tabvalsAMpom}, respectively Table~\ref{tabvalsAMpoa}.  
Here the number of spatial grid points has been chosen such that the smallest spatial mesh width is about equal to $0.40$ and the MCS2-IT$(2)$ 
method has been applied with time step size $\Delta t = 0.01$.
We find a good agreement with the approximations for the American put-on-the-min values in the case of Set~1 obtained in \cite[p.766]{CF08}.
The maximal absolute error in the approximations given by Tables~\ref{tabvalsAMpom} and~\ref{tabvalsAMpoa} is estimated to be less than 0.01. 

\section{Conclusion}\label{sec:conc}
We have investigated a variety of temporal discretization methods in the numerical solution of the two-dimensional time-dependent 
partial integro-differential complementarity problem for the values of American-style options under the two-asset Merton jump-diffusion 
model.
The methods under consideration are constructed as adaptations of modern operator splitting schemes of the IMEX and ADI kind for 
partial integro-differential equations, relevant to European-style two-asset options.
Here the two-dimensional integral part is always conveniently treated in an explicit fashion.
Their adaptation to American options is achieved by a combination with two popular techniques from the computational finance 
literature, the IT splitting technique and the penalty method.
In the present paper we propose the novel approach where IT splitting is employed in an iterative manner in each time step.
We refer to this as IT$(\kappa)$ splitting, where $\kappa\ge 1$ denotes the number of iterations.

Six IMEX and ADI schemes that have recently been studied for European two-asset options in~\cite{BH19} are applied with IT$(\kappa)$ 
splitting: CNFI, IETR, CNAB, MCS, MCS2 and SC2A.
In addition, the CNFI and MCS schemes are combined with the penalty method, as introduced in \cite{HFL04} and \cite{HC18}, respectively.

Ample numerical experiments have been performed to investigate the convergence behaviour of the temporal discretization errors of 
the acquired eight methods. 
To render a fair comparison, the number of time steps for each method is selected such that over the whole time interval 
$[0,T]$ the methods employ the same number of evaluations of the, computationally dominant, integral part.
As a main result we find that taking $\kappa=2$ leads to a regular, monotonic convergence behaviour for the methods based on IT 
splitting, as opposed to the standard application of this technique where $\kappa=1$.
The obtained numerical orders of convergence with $\kappa=2$ lie in between 1.2 and 1.9.

The CNFI-P method leads to temporal errors that are similar to those obtained with all IT(2) splitting methods together, but may 
show a less regular convergence behaviour.
The MCS-P method often yields relatively large temporal errors.

In view of the above and taking into account the size of the error constants, the CNAB-IT$(2)$ and MCS2-IT$(2)$ methods merit 
our preference for the efficient and stable temporal discretization of the two-dimensional Merton PIDCP.

\section*{Acknowledgements}
The authors acknowledge the support of the Research Fund (BOF) of the University of Antwerp (41/FA070300/3/FFB150337).

\clearpage

\begin{table}[h!]
\centering
\begin{tabular}{ l  c c}
 & put-on-the-min & put-on-the-average\\
 \hline
 CNFI-IT$(2)$ & 1.53 &1.53 \\
 CNAB-IT$(2)$ & 1.22 & 1.34\\
 IETR-IT$(2)$ & 1.60 & 1.58 \\
 MCS-IT$(2)$ & 1.59 & 1.58\\
 MCS2-IT$(2)$ & 1.22 & 1.34\\
 SC2A-IT$(2)$ & 1.17 & 1.32\\
 \hline
 CNFI-P & 1.71 & 1.99\\
 MCS-P & 0.85  & 0.89\\
\end{tabular}
\caption[Order of convergence under parameter set 1]{Convergence orders of all methods, with $\kappa=2$, on the small ROI under parameter Set~1.}\label{OoC1}
\end{table}
\vskip0.3cm
\begin{table}[h!]
\centering
\begin{tabular}{ l  c c}
 & put-on-the-min & put-on-the-average\\
 \hline
 CNFI-IT$(2)$ & 1.87 & 1.89\\
 CNAB-IT$(2)$ & 1.80 &1.81\\
 IETR-IT$(2)$ & 1.86 & 1.83\\
 MCS-IT$(2)$ & 1.87 & 1.83\\
 MCS2-IT$(2)$ & 1.83 & 1.82\\
 SC2A-IT$(2)$ & 1.69 & 1.85\\
 \hline
 CNFI-P & 2.04 & 1.59\\
 MCS-P & 1.23 &1.04 \\
\end{tabular}
\caption[Order of convergence under parameter set 2]{Convergence orders of all methods, with $\kappa=2$, on the small ROI under parameter Set~2.}\label{OoC2}
\end{table}
\vskip0.3cm
\begin{table}[h!]
\centering
\begin{tabular}{ l  c c}
 & put-on-the-min & put-on-the-average\\
 \hline
 CNFI-IT$(2)$ & 1.88 & 1.86\\
 CNAB-IT$(2)$ & 1.91 & 1.91\\
 IETR-IT$(2)$ & 1.87 & 1.85\\
 MCS-IT$(2)$ & 1.85 & 1.84\\
 MCS2-IT$(2)$ & 1.93 &  1.93\\
 SC2A-IT$(2)$ & 1.93& 1.93\\
 \hline
 CNFI-P & 1.98 &  1.97\\
 MCS-P & 1.85 & 1.64\\
\end{tabular}
\caption[Order of convergence under parameter set 3]{Convergence orders of all methods, with $\kappa=2$, on the small ROI under parameter Set~3.}\label{OoC3}
\end{table}

\clearpage
\begin{table}[h!]
\centering
\begin{adjustbox}{width = 0.5\textwidth}
\begin{tabular}{ l  r   r  r }
\multicolumn{4}{c}{Set 1} \\
\hline
 & $S_0^{(1)} = 90$ & $S_0^{(1)} = 100$ & $S_0^{(1)} = 110$ \\[0.5cm]
 
$S_0^{(2)} = 90$ & 16.391 & 13.999 & 12.758 \\
$S_0^{(2)} = 100$ & 13.021 & 9.620 & 7.877  \\
$S_0^{(2)} = 110$ & 11.443 & 7.227 & 5.132 \\
\hline
 \end{tabular}
 \end{adjustbox}
 \hspace{0.1cm}\\[0.5cm]
 \begin{adjustbox}{width = 0.5\textwidth}
 \begin{tabular}{ r  r   r  r   }
\multicolumn{4}{c}{Set 2} \\
\hline
 & $S_0^{(1)} = 36$ & $S_0^{(1)} = 40$ & $S_0^{(1)} = 44$\\[0.5cm]
 
$S_0^{(2)} = 36$ & 15.467 & 14.564 & 13.794 \\
$S_0^{(2)} = 40$ & 14.092 & 13.107 & 12.263 \\
 $S_0^{(2)} =44$ & 12.921 & 11.877 & 10.982 \\
\hline
 \end{tabular}
  \end{adjustbox}
   \hspace{0.1cm}\\[0.5cm]
   \begin{adjustbox}{width = 0.5\textwidth}
 \begin{tabular}{ r  r   r  r   }
\multicolumn{4}{c}{Set 3} \\
\hline
 & $S_0^{(1)} = 36$ & $S_0^{(1)} = 40$ & $S_0^{(1)} = 44$\\[0.5cm]
 
$S_0^{(2)} = 36$ & 21.742 & 20.908 & 20.167 \\
$S_0^{(2)} = 40$ & 21.272 & 20.394 & 19.611 \\
 $S_0^{(2)} =44$ & 20.892 & 19.983 & 19.166 \\
\hline
 \end{tabular}
  \end{adjustbox}
\caption{American put-on-the-min option value approximations under parameter Sets~1, 2, 3.}\label{tabvalsAMpom}
\end{table}

\clearpage
\begin{table}[h!]
\centering
\begin{adjustbox}{width = 0.5\textwidth}
\begin{tabular}{ l  r   r  r }
\multicolumn{4}{c}{Set 1} \\
\hline
 & $S_0^{(1)} = 90$ & $S_0^{(1)} = 100$ & $S_0^{(1)} = 110$ \\[0.5cm]
 
$S_0^{(2)} = 90$ & 10.003 & 5.989 & 3.441 \\
$S_0^{(2)} = 100$ & 6.030 & 3.442 & 1.887  \\
$S_0^{(2)} = 110$ & 3.491 & 1.891 & 0.993 \\
\hline
 \end{tabular}
 \end{adjustbox}
 \hspace{0.1cm}\\[0.5cm]
 \begin{adjustbox}{width = 0.5\textwidth}
 \begin{tabular}{ r  r   r  r   }
\multicolumn{4}{c}{Set 2} \\
\hline
 & $S_0^{(1)} = 36$ & $S_0^{(1)} = 40$ & $S_0^{(1)} = 44$\\[0.5cm]
 
$S_0^{(2)} = 36$ & 5.406 & 4.363 & 3.547 \\
$S_0^{(2)} = 40$ & 4.214 & 3.339 & 2.669 \\
 $S_0^{(2)} =44$ & 3.225 & 2.507 & 1.969 \\
\hline
 \end{tabular}
  \end{adjustbox}
   \hspace{0.1cm}\\[0.5cm]
   \begin{adjustbox}{width = 0.5\textwidth}
 \begin{tabular}{ r  r   r  r   }
\multicolumn{4}{c}{Set 3} \\
\hline
 & $S_0^{(1)} = 36$ & $S_0^{(1)} = 40$ & $S_0^{(1)} = 44$\\[0.5cm]
 
$S_0^{(2)} = 36$ & 12.466 & 11.930 & 11.440 \\
$S_0^{(2)} = 40$ & 11.434 & 10.943 & 10.495 \\
 $S_0^{(2)} =44$ & 10.493 & 10.043 & 9.633 \\
\hline
 \end{tabular}
  \end{adjustbox}
\caption{American put-on-the-average option value approximations under parameter Sets~1, 2, 3.}\label{tabvalsAMpoa}
\end{table}

\clearpage
\bibliographystyle{plain}
\bibliography{bib2DMerton_splitting}

\end{document}